\font\script=eusm10.
\font\sets=msbm10.
\font\stampatello=cmcsc10.
\def\0{{\bf 0}}
\def\1{{\bf 1}}
\def\defineq{\buildrel{def}\over{=}}
\def\definiz{\buildrel{def}\over{\Longleftrightarrow}}

\def\C{\hbox{\sets C}}
\def\N{\hbox{\sets N}}
\def\R{\hbox{\sets R}}
\def\P{\hbox{\sets P}}

\def\cloudF{\hbox{$<F>$}}
\def\cloud0{\hbox{$<\0>$}}
\def\corsivoC{\hbox{\script C}}
\def\corsivoD{\hbox{\script D}}
\def\corsivoE{\hbox{\script E}}
\def\corsivoF{\hbox{\script F}}
\def\corsivoH{\hbox{\script H}}
\def\corsivoL{\hbox{\script L}}
\def\corsivoM{\hbox{\script M}}
\def\corsivoR{\hbox{\script R}}
\def\corsivoS{\hbox{\script S}}
\def\corsivoT{\hbox{\script T}}
\def\corsivoU{\hbox{\script U}}

\def\square{\hbox{\vrule\vbox{\hrule\phantom{s}\hrule}\vrule}}

%
\def\qed{\hfill $\square$\par}

%
%
%
\def\CGf{{\corsivoC_G}}  
\def\DGf{{\corsivoD_G}}
\def\EGf{{\corsivoE_G}}
\def\LGf{{\corsivoL_G}}
\def\MGf{{\corsivoM_G}}
\def\RGf{{\corsivoR_G}}
\def\SGf{{\corsivoS_G}}
\def\UGf{{\corsivoU_G}}
\def\CG#1{{\corsivoC_G(#1)}}  
\def\DG#1{{\corsivoD_G(#1)}}
\def\EG#1{{\corsivoE_G(#1)}}
\def\LG#1{{\corsivoL_G(#1)}}
\def\MG#1{{\corsivoM_G(#1)}}
\def\RG#1{{\corsivoR_G(#1)}}
\def\ETRG#1{{\corsivoR'_G(#1)}}
\def\SG#1{{\corsivoS_G(#1)}}
\def\UG#1{{\corsivoU_G(#1)}}
\def\RGx#1#2{{\corsivoR_G(#1,#2)}} 
\def\LGx#1#2{{\corsivoL_G(#1,#2)}}
\def\SGx#1#2{{\corsivoS_G(#1,#2)}}
%

%
\def\corsivoA{\hbox{\script A}}
\def\Bad{\hbox{\script B}}
\def\simplyB{{\Bad^s}}
\def\hyperB{{\Bad^\infty}}
\def\transparent{\corsivoT}
\def\simplyT{{\transparent^s}}
\def\hyperT{{\transparent^\infty}}
\def\valg{{v_{p,G}}}
\def\wg{{w_{p,G}}}

\def\rad{\mathop {\rm rad}}

\def\limx{\lim_{x\to\infty}}
%
%

\def\sumq{\sum_{q=1}^\infty}
\def\sumr#1{\sum_{(r,#1)=1}}

\def\and{\quad\quad\hbox{\stampatello and}\quad\quad}
\def\nondivide{\!\not\,\mid }	
\def\Hildebrand{{\rm Hi}_F}	
%
\def\arithf{\,:\, \N \rightarrow \C}
\def\cloudF_M{<F>_M}
\def\gmf{G_{M_F}}
\def\nt{N_T(G)}
\def\SH{\corsivoS_H}		
\def\vpa{v_p(a_F)}
\par
\centerline{\bf Convergence of Ramanujan expansions, I}
\medskip
\centerline{\rm [Multiplicativity on Ramanujan clouds]}
\bigskip
\medskip
\bigskip
\centerline{Giovanni Coppola and Luca Ghidelli}
\bigskip
\par
\noindent
{\bf Abstract}. 
We call \enspace $\RG a:=\sum_{q=1}^{\infty}G(q)c_q(a)$ \enspace the {\stampatello Ramanujan series}, of coefficient $G:\N\to\C$, where $c_q(a)$ is the well-known {\it Ramanujan sum}. We study the convergence of this series (a preliminary step, to study {\it Ramanujan expansions}) and define $G$ a {\it Ramanujan coefficient} when $\RG a$ converges pointwise, in all natural $a$. Then, $\RGf:\N\to\C$ is {\it well defined} ({\stampatello w-d}). 
The {\stampatello Ramanujan cloud} of a fixed $F:\N\to\C$ is $<F>:=\{G\!:\!\N\!\to\!\C \;|\; \RGf \enspace \hbox{\stampatello w-d}, F=\RGf\}$. (See the Appendix.) We study in detail the multiplicative Ramanujan coefficients $G$ : their $<F>$ subset is called the {\stampatello multiplicative Ramanujan cloud}, $<F>_M$. 
\par
Our first main result, the \lq \lq {\stampatello finiteness convergence Theorem}\rq \rq, for $G$ multiplicative, among other properties equivalent to \lq \lq $\RGf$ well defined\rq \rq, reduces the convergence test to a finite set, i.e., $\RGf$ {\stampatello w-d} is equivalent to: $\RG a$ converges for all $a$ dividing $N(G)\in\N$, that we call the \lq \lq {\stampatello Ramanujan conductor}\rq \rq. 
\par
Our second main result, the \lq \lq {\stampatello finite Euler product explicit formula}\rq \rq, for multiplicative Ramanujan coefficients $G$, writes $F=\RGf$ as a finite Euler product; thus, $F$ is a semi-multiplicative function (following Rearick definition) and this product is the Selberg factorization for $F$. In particular, we have: $F(a)=\RGf(a)$ converges absolutely, being {\it finite} (of length depending on non-zero $p-$adic valuations of $a$). 
\par
Our third main result, the \lq \lq {\stampatello multiplicative Ramanujan clouds}\rq \rq, studies the important subsets of $<F>_M$; also giving, for all  multiplicative $F$, the {\it canonical Ramanujan coefficient} $G_F\in <F>_M$, proving: {\stampatello any multiplicative} $F$ {\stampatello has a finite Ramanujan expansion with multiplicative coefficients}. 

\bigskip

\par
\noindent{\bf Table of Contents}
\bigskip

\def\tablewidth{\hsize}
\tabskip=1em  
\halign to \tablewidth 
{
\hfil# \tabskip=1em plus 1em  &
#\hfil \tabskip=2em plus 2em &
\hfill#\cr
1. & Introduction and main results & \bf 1 \cr
2. & Notation, definitions \& basic formul\ae & \bf 6 \cr
3. & Ramanujan, Lucht and coprime series: from infinite to finite Euler products & \bf 9 \cr
4. & Eratosthenes transforms and completely multiplicative clouds & \bf 14 \cr
5. & Canonical Ramanujan coefficients and the separating multiplicative clouds & \bf 15 \cr
6. & Converse convergence Theorem & \bf 19 \cr
7. & Finiteness convergence Theorem: Proof of Theorem 1.1 & \bf 21 \cr
8. & Ramanujan series factorization and Euler-Selberg products & \bf 22 \cr
9. & Multiplicative Ramanujan clouds & \bf 27 \cr
Appendix. & Odds \& ends (more details, alternative proofs, generalizations) & \bf 33 \cr
}

\bigskip
\bigskip
\bigskip

\par
\noindent{\bf 1. Introduction and Main results}
\bigskip
\par
\noindent
We began studying the Ramanujan expansions [R] from a new point of view. In fact, usually, we have a fixed arithmetic function $F:\N\to\C$ (as usual, $\N$ and $\C$ are natural and complex numbers) and we wish to, if possible, find a sequence of \lq \lq Ramanujan coefficients\rq \rq, $G(1),G(2),\ldots,G(q),\ldots$, say {\it a Ramanujan coefficient} $G:\N\to \C$ such that 
$$
F(n)=\sum_{q=1}^{\infty}G(q)c_q(n), 
$$
\par
\noindent
say, for all the natural numbers $n$. Here, as usual, the trigonometric sum $c_q(n)$ is the {\it Ramanujan sum} [R] : 
$$
c_q(n)\defineq \sum_{j\le q, (j,q)=1}\cos{{2\pi jn}\over q} , 
$$
\par
\noindent
where we write $(a,b)$ as an abbreviation for their g.c.d. : so, $(j,q)=1$ means they're coprime. 
\bigskip
\par				
First of all, we know that, for each fixed $F:\N\to\C$, there are a lot (uncountably many!) Ramanujan coefficients, compare [C3]. 
\par
Second, the problem of finding, once fixed $F$, {\bf all} $G$ such that the above expansion holds, needs of course, first, to {\bf check the convergence} ! However even this (seemingly) easier problem, in absence of properties for $G$, could be a hopeless, at least huge, task. So, in the present paper we will confine ourselves with {\it multiplicative Ramanujan coefficients} $G$. (This general rule will be relaxed in the Appendix.) 
\par
In the following, we'll abbreviate \lq \lq if and only if\rq \rq, with \lq \lq iff\rq \rq.
\par
\noindent
We define $G$ {\it multiplicative} iff $G(ab)=G(a)G(b)$, for all coprime $a,b\in\N$; so, including the constant $0-$function, or {\it null-function}, say $\0(a)\defineq 0$, $\forall a\in\N$ (however, $G\neq \0$ multiplicative $\Rightarrow$ $G(1)=1$). 
\par
In future papers, but with a glance already in this (compare the Appendix), we will face the problem of finding, once fixed $F:\N\to\C$, ALL its Ramanujan coefficients $G$. 
\par
In all, our approach will be to fix a $G:\N\to\C$, in this paper (except Appendix) multiplicative, and first check for the convergence in $\N$ of the corresponding {\stampatello Ramanujan series} (compare [C3])
$$
\RG a \defineq \sum_{q=1}^{\infty}G(q)c_q(a), 
\enspace 
\forall a\in \N; 
$$
\par
\noindent
then, try to describe, from the hypotheses on $G$, the properties of its sum, say $F(a)=\RG a$. 
\par
An important warning : hereafter, when we say \lq \lq {\it Ramanujan coefficient}\rq \rq, we mean any arithmetic function $G:\N\to\C$, such that Ramanujan series $\RG a$ {\bf converges pointwise in all natural numbers} $a$. (For example, in [La] the Author considers also convergence in subsets of $\N$.)
\par
In other words, we wish to study the so-called {\stampatello Ramanujan cloud} [C3] of a fixed $F:\N\to\C$, that {\it is the set of its Ramanujan coefficients}$\defineq <F>$; but, since this is too difficult (at least nowadays), we study {\it parts of a cloud}, namely we look for subsets, of this $<F>$. 
\par
Let us say, the most beautiful subset of a cloud (up to now, also considering [C3] subsets, like $<F>_{\ast}$, $<F>_{\ast \ast}$ and $<F>_{\#}$, compare the Appendix) is the so-called {\stampatello Multiplicative Ramanujan Cloud}, i.e.: 
$$
<F>_M\defineq \{G\in <F>\; :\; G\enspace \hbox{\stampatello multiplicative}\}.  
$$
\par
\noindent
Our main concern is the study of this $<F>_M$ and, also, of how the $F$ properties and the $G$ properties influence each other. One typical question, for example: what properties, of a fixed $F$, ensure $<F>_M\neq \emptyset$ ? 
\medskip
\par
\noindent
The first thing to check is the convergence of $\RGf$ and so our first result gives a set of properties, all equivalent to: $G$ is a Ramanujan coefficient (recall, one last time, meaning : $\RG a$ converges $\forall a\in\N$). 
\par
A kind of first elementary approach to this problem began in [C2], in which the first Author studied the \lq \lq finite Euler products\rq \rq, of Ramanujan series. There, we encounter the {\it coprime series}, with the usual {\it M\"obius function} $\mu$ (see $\S2$)
$$
\SG a \defineq \sum_{(r,a)=1}G(r)\mu(r),
\enspace 
\forall a\in \N.
$$
\medskip
\par
Then, the Authors started to study the multiplicative cloud of $\0$, the null-function (see our [CG]), but under hypotheses of convergence, for these coprime series. Here, thanks to a wider view on convergence conditions (compare next Theorem 1.1), we are able to complete this $<\0>_M$  Classification, see $\S9.1$. 
\medskip
\par
Then, we realized that a much better understanding of Ramanujan expansions convergence might come from classic general results (of pointwise convergence), due mainly to Lucht [Lu1], compare the survey [Lu2]. In building upon Lucht's insight, compare Corollary 3.7 and see $\S4$, we define a new series : the {\it Lucht series} 
$$
\LG d \defineq \sum_{K=1}^{\infty}\mu(K)G(dK),
\enspace 
\forall d\in \N. 
$$
\par
\noindent
Notice an important link, namely (in case of common convergence in $1$) 
$$
\RG 1 = \LG 1 = \SG 1. 
$$

\bigskip

\par				
We are now ready to state our main technical theorem concerning the convergence of Ramanujan expansions with multiplicative coefficients. 
Notice that the necessary condition of convergence of $\RG 1 = \SG 1 = \LG 1$ implies (esp., compare Remark 1 in [C2] and, also, see [CG]) having finitely many bad primes, see Definition A in $\S2$; there, the Ramanujan conductor $N(G)$ and hyperbad primes, $p\in\hyperB$, are defined. 
\smallskip
\par
\noindent {\bf Theorem 1.1.} ({\stampatello Finiteness convergence theorem})
{\it Let \thinspace $G$ \thinspace be a multiplicative function, such that \enspace $\RG 1 = \SG 1 = \LG 1$ \enspace converges. 
Then the following are equivalent:
\smallskip
\quad (1) \enspace $\RG a$ converges for all $a\in\N$;
\smallskip
\quad (2) \enspace $\RG a$ converges for all $a\mid N(G)$;
\smallskip
\quad (3) \enspace $\SG {N(G)}$ converges;
\smallskip
\quad (4) \enspace $\SG b$ converges for all $b$ not divisible by any hyperbad prime $p\in\hyperB$. 
} 
\bigskip
\par
\noindent 
We prove this theorem in Section 7, via the sequence of implications : (1)$\Rightarrow$(2)$\Rightarrow$(3)$\Rightarrow$(4)$\Rightarrow$(1). 
\par
\noindent 
The implication (1)$\Rightarrow$(2) is trivial, while the implication (3)$\Rightarrow$(4) will follow from the \lq \lq Converse convergence Theorem\rq \rq, that we prove in Section 6. 
\par 
The formul\ae, for Ramanujan series $\RGf$ and coprime series $\SGf$, provided in Section 3 will prove, in particular from Corollaries 3.3 and 3.4, the implication (4)$\Rightarrow$(1). 
\par
The implication (2)$\Rightarrow$(3) is the \lq \lq most difficult\rq \rq \enspace part of the proof; in fact, it is based on a formula (see $\S7$) at the heart of the links in between $\RGf$, $\SGf$ and $\LGf$, involving the three variables function $\corsivoF_G$, see $\S2$. This key formula is stated in Theorem 7.1. 
\par
The complete proof of Theorem 1.1 is given in $\S7$. 
\bigskip
\par
\noindent
After we know how to characterize the convergence in $\N$, with previous result, we wish to know which kind of function is $\RGf:\N\to\C$ (always well-defined in $\N$ hereafter). 
\par
For this, we give Euler products formul\ae, that might seem expansions into infinite products, but at a closer inspection reveal to be, actually, {\it finite Euler products} ! A kind of unexpected result is that multiplicative Ramanujan coefficients $G$ give Ramanujan expansions $\RGf$ that are not multiplicative functions, but only slightly more general arithmetic functions: the $\RGf$ {\it are semi-multiplicative}. (Compare section end, where semi-multiplicative functions are defined and also characterized in terms of finite Euler products.) 
\medskip
\par
The \lq \lq Ramanujan series factorization\rq \rq, Theorem 8.1 in $\S8$, implies in fact a kind of explicit form in finite Euler products. Recall, from $\S2$, the $p-$adic valuations $v_p(a)$ and $v_{p,G}$, like Definitions A and B. 
\smallskip
\par
\noindent {\bf Corollary 1.2.} ({\stampatello finite Euler products explicit formula})
\par
\noindent 
{\it Let $G:\N\to\C$ be a multiplicative function such that $\RGf$ converges in $\N$. Then $\RG a = 0$ for every $a\not \equiv 0\,\bmod N_T(G)$. Moreover, there is a multiplicative function $\MGf : \N\to\C$ such that for all $a\in\N$: 
$$
\RG {aN_T(G)} = \RG {N_T(G)}\cdot \MG a
= \RG {N_T(G)}\cdot \prod_{p|a} \MG {p^{v_p(a)}}. 
$$
\par
\noindent
Hence, we have trivially : 
$$
\RG {N_T(G)} = 0
\enspace \Rightarrow\enspace 
\RGf = \0.
$$ 
\par
\noindent 
If \enspace $\RG {N_T(G)}\neq 0$, then the function $\MGf$ is uniquely determined, with the following values at prime powers: 
\medskip
\quad $\MG {p^K} \defineq 1 + p + \dots + p^K$, \enspace {\it if $p$ is hypertransparent}; 
\medskip
\quad $\MG {p^K} \defineq (1-G(p^{\valg+1}))^{-1} {\displaystyle \sum_{k=0}^K  p^{k} (G(p^{\valg+k}) - G(p^{\valg+k+1})) }$, \enspace {\it if $p$ is simply transparent};
\smallskip
\quad $\MG {p^K} \defineq (1-G(p))^{-1}\; {\displaystyle \sum_{k=0}^{K}p^k\left(G(p^k)-G(p^{k+1})\right) }$, \enspace {\it otherwise}. 
\par
\noindent 
Furthermore, in the hypothesis $G$ multiplicative Ramanujan coefficient, we get the general formula: 
$$
\RG a 
= \RG {N_T(G)} 
  \cdot 
\prod_{{\valg=\infty}} \left( \sum_{v=0}^{v_p(a)} p^v \right)
   \cdot 
\prod_{{\valg\neq \infty}} \left(\sum_{v=0}^{v_p(a)} p^{v-\valg} { G(p^v)- G(p^{v+1}) \over 1 - G(p^{\valg +1}) } \right), 
\quad \forall a\in \N. 
$$
\par				
\noindent 
Note that all the above products are finite: all but finitely many factors are equal to $1$, because $\valg=v_p(a)=0$ for all but finitely many primes. 
Also, when $\RG {N_T(G)} = 0$, the value of this right-hand side is unimportant, because of previous trivial implication. 
} 
\medskip
\par
\noindent
We prove this Corollary in $\S8$, after having proved its, say, relative Theorem 8.1(in section 8). 

\bigskip

\par
We give here the necessary definitions and properties of semi-multiplicative and Selberg-multiplicative arithmetic functions, that we apply to interpret our results. Recall, from $\S2$, that $\N_0=\N\cup \{0\}$. 
\smallskip
\par
The {\it semi-multiplicative} arithmetic functions $F:\N\to\C$ where defined by Rearick [Re] as those satisfying, for a certain complex $c\neq 0$ and an integer $a_F>0$, 
$$
F(n)=cM_F(n/a_F)
\enspace \forall n\in \N, 
$$
\par
\noindent
where the function $M_F$, vanishing outside $\N$ (so $F(n)=0$, on all $n\not\equiv 0\bmod a_F$), is multiplicative. Notice that $F=\0$ iff $M_F=\0$ and that, whenever $F\neq \0$, we have $c=F(a_F)$, where $a_F=\min\{v\in\N : F(v)\neq 0\}$ is called the {\it threshold}. Of course, $F\neq \0$ multiplicative implies $F$ semi-multiplicative, with $a_F=1$ and $c=1$. 
\par
Another generalization of \lq \lq $F:\N\to\C$ multiplicative\rq \rq \thinspace is when $F$ is {\it Selberg-multiplicative}: 
$$
F(a)=\prod_{p\in \P}F_p(v_p(a)), 
\forall a\in \N 
$$
\par
\noindent
(recall : $p-$adic valuation $v_p(a)$ in $\S2$), where the functions $F_p:\N_0\to\C$ have $F_p(0)=1$ for all but finitely many primes $p$. The importance of this {\it Selberg factorization} comes from the explicit form we can give to its factors $F_p$ and, also, that Selberg multiplicativity is the same as semi-multiplicativity: both these properties were proved, by Haukkanen, see [Hau], Theorem 2.8, that we quote, here. 
\medskip
\par
\noindent {\bf Proposition 1.3} ({\stampatello Semi-multiplicative $F$ are Selberg-multiplicative : explicit factorization})
\par
\noindent
{\it An arithmetic function $F:\N\to\C$ is Selberg multiplicative iff it is semi-multiplicative. In this case,
$$
F(a)=F(a_F)\prod_{p\in\P}\left({{F(a_F p^{v_p(a)-v_p(a_F)})}\over {F(a_F)}}\right)
$$
\par
\noindent
is the Selberg factorization of our $F$.
}

\bigskip

\par
\noindent
In fact : take $F=\RGf$, $a_F=N_T(G)\ge 1$ and $c=F(a_F)=\RG {N_T(G)}\neq 0$ (otherwise, $F=\0$) in Corollary 1.2; so, it proves, in particular, that, for multiplicative Ramanujan coefficients $G$, the Ramanujan expansion $\RGf$ is a semi-multiplicative function, with above Selberg factorization. Notice, once again, that it is, of course, a {\it finite} Euler product. Sometimes we'll call it the {\stampatello Euler-Selberg product} of our Ramanujan expansion. 

\bigskip

\par
We arrive to the structure of multiplicative Ramanujan clouds. Next Corollary is a kind of summary of all our results on multiplicative Ramanujan clouds (see $\S9$) and its important subsets ($\S4,\S5$), for which see $\S2$ definitions. Thus we will not prove it. 
\smallskip
\medskip
\par
\noindent {\bf Corollary 1.4} ({\stampatello Multiplicative Ramanujan Clouds})
\par
\noindent
{\it We start with the easiest multiplicative Ramanujan cloud: that of\enspace $\0$, the null-function: 
$$
<\0>_M = \{ G:\N\to\C \;|\; G\enspace \hbox{\stampatello multiplicative,\enspace with }\SG {N(G)} = 0\}. 
$$
\par
\noindent
Let $F\neq \0$ hereafter. 
\par
Then $G$ is multiplicative $\Longleftrightarrow$ \enspace $\RGf=F$ \thinspace is semi-multiplicative. (The $\Rightarrow$ from Corollary 1.2 and the $\Leftarrow$ from Proposition 9.8.). In this case, $\RG {N_T(G)} \neq 0$; otherwise, by previous Corollary 1.2, $F=\0$. 
\par
\noindent
Coming to remarkable subsets of $<F>_M$ (recall their definitions in $\S2$), consider first the Euler separating multiplicative Ramanujan cloud of a fixed, non-null multiplicative $F:\N\to\C$, i.e. 
$$
F\neq \0 \enspace \hbox{\stampatello is\enspace multiplicative}
\quad \Longleftrightarrow \quad
<F>_{ESM}=\{G_F\}\,, 
$$
\par				
\noindent
where the $G_F$ is the {\stampatello canonical Ramanujan coefficient} of this $F$, see $\S2$. In particular, (compare $\S A.6$)
\smallskip
\par
\centerline{\stampatello $F\neq \0$ is multiplicative \thinspace $\Longleftrightarrow $ \thinspace $F=\RGf$ is finite and pure with multiplicative $G=G_F\,.$}
\smallskip
\par
\noindent
The separating multiplicative Ramanujan cloud of a fixed $F:\N\to\C$ is, as follows, 
$$
<F>_{SM} \neq \emptyset
\quad \Longleftrightarrow \quad
F\enspace \hbox{\stampatello is\enspace quasi-multiplicative}. 
$$
\par
\noindent
In fact, the multiplicative Ramanujan coefficients in it, say $G$, have $\SGf/F(1)$, $F/F(1)=\RGf/F(1)$ and $\LGf/F(1)$ that are all well defined (with $\SG 1=\RG 1=\LG 1=F(1)\neq 0$, of course) and multiplicative. 
\par
\noindent
Thus, in particular, 
$$
<F>_{NSM}\neq \emptyset
\quad \Longleftrightarrow \quad 
F\enspace \hbox{\stampatello is\enspace multiplicative}. 
$$
\par
\noindent
The Ramanujan coefficients, say $G$, in the multiplicative Ramanujan cloud of a semi-multiplicative $F$, i.e., \enspace $G\in <F>_M$, \enspace are determined (inside $<F>_M$) by a square-free supported arithmetic function, say $H_G$, called the \lq \lq {\stampatello opacity core}\rq \rq \thinspace of $G$ (see in $\S9.2$,  Theorem 9.3). 
\par
\noindent
Finally, the complete Classification for the $<F>_M$, when $F\neq \0$, is given in Proposition 9.8. 
}

\vfill
\eject

\par				
\noindent{\bf 2. Notation, definitions \& basic formul\ae}
\bigskip
\par
\noindent
Typographic definitions: $\diamond$ is a Remark's end, \lq \lq QED\rq \rq \thinspace is to separate parts of a proof, which ends with \hfill $\square$ 
\par
\noindent
We abbreviate $\N_0\defineq \N \cup\{0\}$, with $\N\defineq $set of {\it natural numbers}; while $\P\defineq$set of {\it primes}, usually denoted $p$. 
\par
An {\it arithmetic function} is any $F:\N\to\C$ (\lq \lq a sequence of complex numbers\rq \rq, outside Number Theory). 
\par
The {\it characteristic function} of a subset $\corsivoH \subseteq \N$ of natural numbers, an arithmetic function, will be as usual: $\1_{\corsivoH}(n)\defineq 1$, iff $n\in \corsivoH$ ($\defineq 0$, otherwise); in particular $\1\defineq \1_{\N}$ is the constant$-1-$function. 
\par
\noindent
The {\it completely multiplicative} arithmetic functions are $f:\N\to\C$ with $f(ab)=f(a)f(b)$, $\forall a,b\in\N$. If this happens when $a,b\in\N$ are coprime, we call it {\it multiplicative}. The {\it null-function} (or $0-${\it constant function}) $\0(a)\defineq 0$, $\forall a\in\N$, is (here) a multiplicative function: the only one with $f(1)=0$, all others have $f(1)=1$. 
\par
\noindent
Also, $\1_{\{1\}}$ will denote the function which is always $0$, except in $1$: $\1_{\{1\}}(1)\defineq 1$, a completely multiplicative arithmetic function. (For $\0,\1,\1_{\{1\}}$, compare [La].) Once fixed any $p$, the classic {\it $p-$adic valuation} is : $v_p(a)\defineq \max\{k\in\N_0 \, :\, p^k|a\}$, $\forall a\in \N$, a {\it completely additive} function: $v_p(ab)=v_p(a)+v_p(b)$, $\forall a,b\in\N$. 
\medskip
\par
\noindent
The {\it Ramanujan sum} of modulus $q\in\N$ and argument $a\in\N$, $c_q(a)$, defined above, has, resp., Kluyver's formula [K], $(3)$ in [M], and H\"older's formula [H\"o], [M] (we'll use both extensively throughout the paper): 
$$
c_q(a)=\sum_{d|q,d|a}d\mu(q/d), 
\quad
c_q(a)=\varphi(q){{\mu(q/(q,a))}\over {\varphi(q/(q,a))}}
\qquad
\forall q,a\in\N. 
$$
\par
\noindent
Here $\mu$ is {\it M\"obius function}, a multiplicative function with $\mu(1)\defineq 1$, $\mu(p)\defineq -1$, $\forall p\in \P$, while $\exists p\in\P$:  $p^2|q$ implies $\mu(q)\defineq 0$. In fact, $\mu=0$ outside the {\it square-free} numbers : these have characteristic function $\mu^2$. 
\par
\noindent
Also, the {\it Euler function} is $\varphi(q)\defineq |\{n\le q\, :\, (n,q)=1\}|$, the cardinality (we'll use $|\corsivoA|\defineq$number of elements, for finite sets $\corsivoA$) of naturals, below $q$, coprime with $q$. Like $\mu$, it is multiplicative. 
\par
\noindent
Given $G:\N\to\C$ multiplicative and $p$ prime, we introduce the quantity $\wg\in\N\cup\{\infty\}$ given by 
$$
\wg\defineq \max \{ w\in\N: \ G(p^k) = G(p)^k, \forall k\le w\},
$$
\par
\noindent
with the convention that $\wg\defineq \infty$ if this set is $\N$, that is, if $G$ is \lq \lq completely multiplicative\rq \rq \thinspace along powers of $p$, see [CG] : we call the $\wg$ defined above the {\it completely multiplicative (C.M.) index} of $p$ for $G$. 
\medskip
We come, now, to all the definitions regarding the primes, with respect to the fixed, multiplicative $G$. 
\medskip
\par 
\noindent {\bf Definition A.} ({\stampatello Bad primes and transparent primes})
\par
\noindent 
{\it Let $G:\N\to\C$ be a given multiplicative function. We say that $p\in\P$ is a {\bf bad prime} for $G$ iff $1\leq |G(p)|\leq p$, and, in particular, a {\bf transparent prime} iff $G(p)=1$. 
We say that $p$ is {\bf simply bad}, resp. {\bf hyperbad}, iff $p$ is bad and $\wg\neq \infty$, resp. is bad and $\wg=\infty$. 
We say that $p$ is {\bf simply transparent}, resp. {\bf hypertransparent}, iff $p$ is transparent and $\wg\neq \infty$, resp. is transparent and $\wg=\infty$. 
The sets of bad, simply bad, hyperbad, transparent, simply transparent, hypertransparent primes for $G$ will be denoted respectively by: $\Bad$, $\simplyB$, $\hyperB$, $\corsivoT$, $\simplyT$, $\hyperT$. All these depend on $G$, but we'll omit $G$ subscript. 
}
\medskip
\par
Recall the {\it transparency index} of $p$ for $G\neq \0$ that we defined in [CG] 
$$
\valg \defineq \min\{ K\in\N_0 \colon \  G(p^{K+1})\neq 1\},
$$
\par
\noindent
with the convention that $\valg=\infty$ if this set is empty (that is, if $p$ is hypertransparent). 
In general, the quantity $\valg$ should not be confused with $\wg$: if $G(p)\neq 1$ we have $\valg=0$, whereas $\wg\geq 1$ for all $p$ prime. 
Note however that $\valg=\wg$ if $G(p)=1$ (that is, when $p$ is a transparent prime). Recall from [CG]: an arithmetic function $G$ is {\stampatello  normal} iff (by definition) it's multiplicative and has no transparent primes. Last but not least, see that we sometimes call $\valg$ the {\it $p-$adic valuation of $G$} (on transparent primes). 

\bigskip

\par
Next, we introduce two notions of \lq \lq conductors\rq \rq, which encode the information coming from simply bad and simply transparent primes. 
\smallskip
\par				
\noindent {\bf Definition B.} ({\stampatello Ramanujan conductor and transparency conductor})
\par
\noindent
{\it If $G$ is a multiplicative function with finitely many bad primes, we define its {\bf Ramanujan conductor} $N(G)$ and its {\bf transparency conductor} $N_T(G)$ via the C.M. index \enspace $\wg\enspace :$ 
$$
N(G)\defineq \prod _{p\in\simplyB} p^{\wg}
\and 
N_T(G)\defineq \prod _{p\in\simplyT} p^{\wg}
=\prod _{p\in\simplyT} p^{\valg}.
$$ 
Note that the products above are restricted to simply bad and simply transparent primes respectively. 
If $G$ has infinitely many bad primes, we agree that \enspace $N(G)\defineq N_T(G)\defineq 0$. 
}

\bigskip

We recall a classical object, in Number Theory, whose name is due to Wintner. 
\smallskip
\par
\noindent {\bf Definition C.} ({\stampatello Eratosthenes Transforms of arbitrary arithmetic functions})
\par
\noindent
{\it Given any arithmetic function $f:\N\to\C$, its } {\stampatello Eratosthenes Transform} ({\it Wintner's $[Wi]$ terminology}) {\it is 
$$
f'\defineq f\ast \mu, 
$$
\par
\noindent
where $\ast$ is the usual {\stampatello Dirichlet product} $(see\enspace [T])$ between arithmetic functions; by M\"obius inversion $[T]$, this is equivalent to:
$$
f=f'\ast \1,
\enspace \hbox{\rm i.e.}, \enspace
f(n)=\sum_{d|n}f'(d),
\forall n\in\N. 
$$
}
\par
\noindent
Recall: $f$ is multiplicative iff $f'$ is multiplicative [T]. \enspace By the way, every {\it divisor} is a natural number, here. 
\medskip
\par
A set (of natural numbers) may be very useful, when passing from some naturals to their divisors. 
\smallskip
\par
\noindent {\bf Definition D.} ({\stampatello divisor-closed subsets of natural numbers})
\par
\noindent
{\it We call } $\corsivoD\subseteq \N$, $\corsivoD\neq \emptyset$ {\it a divisor-closed set iff } 
$$
a\in \corsivoD, d|a
\enspace \Longrightarrow \enspace
d\in \corsivoD. 
$$

\bigskip

\par
\noindent
Throughout the paper it will be useful to use truncated versions of the above series (we'll also drop, at some point, the subscript \lq \lq $G$\rq \rq, when the multiplicative function $G$ is implicit): we define, $\forall a\in \N$, $\forall x\ge 0$, 
$$
\RGx a x \defineq \sum_{q\leq x} G(q) c_q(a); 
\enspace \thinspace \enspace 
\LGx a x \defineq \sum_{K\leq x}\mu(K)G(aK); 
\enspace \thinspace \enspace 
\SGx a x \defineq \sum_{{r\leq x \atop (r,a)=1}}G(r)\mu(r). 
$$
\par
\noindent
(Of course, $0\le x<1$ \enspace gives empty sums, which are $0$ : however, we consider all real $x\ge 0$, for technical reasons in the Converse Convergence Theorem $\S6$. Also, this convention simplifies details in formul\ae.) 
\par
A kind of unification of these partial sums is, for $G:\N\to\C$ and the triple $(a,b,c)\in\N^3$ both fixed: 
$$
\corsivoF_G(a,b,c)(x)\defineq \sum_{{q\le x}\atop {(q,b)=1}}G(cq)c_q(a), 
$$
\par
\noindent
which, if in the limit $x\to \infty$ converges (in $\C$, for the same $(a,b,c)\in \N^3$ and $G$ fixed), gives the \lq \lq $\corsivoR \corsivoS \corsivoL$ series\rq \rq: 
$$
\corsivoF_G(a,b,c)\defineq \sum_{{q=1}\atop {(q,b)=1}}^{\infty}G(cq)c_q(a). 
$$
\par
\noindent
Notice : this case (i.e., of convergence) entails that all of our three series,
$$
\corsivoF_G(a,1,1)=\RG a,
\quad
\corsivoF_G(1,b,1)=\SG b,
\quad
\corsivoF_G(1,1,c)=\LG c,
$$
\par				
\noindent
converge. Even with respect to the two series in [La], this $\corsivoF_G$ is more general : $\Phi_G(a,b)=\corsivoF_G(a,b,1)$ and, since $\Phi_G(a,b)=\sum_{d|b}\mu(d)\Phi_G'(a,d)$, by M\"obius inversion [T], $\Phi_G'(a,d)=\mu(d)\sum_{b|d}\corsivoF_G(a,b,1)\mu(d/b)$, for all square-free $d\in\N$. 
\par
We'll prove the most difficult part of previous Theorem 1.1 by a transformation formula for $\corsivoF_G$, see $\S7$. 
\medskip
\par
\noindent
We pass now to the definition of finite sums, playing the part of finite factors (in Euler products formul\ae). 
\medskip
We define two sums : $\corsivoE_G(a)$, the {\it Euler-Ramanujan factor}, and, resp., $\corsivoU_G(a)$, the {\it Lucht-Ramanujan factor}, that appear in the factorization of $\RG a$ and, resp., $\LG a$. For any multiplicative $G:\N \rightarrow \C$, given any $a\in \N$, $\forall p\in \P$, recall $v_p(a)$ is the $p-$adic valuation of $a$; {\it the $p-$Euler factor} \thinspace of \thinspace $G(q)c_q(a)$ \thinspace is: 
$$
E_{p,G}(a)\defineq \sum_{K=0}^{v_p(a)+1}G(p^K)c_{p^K}(a)
\enspace \Rightarrow \enspace 
\EG a \defineq \sum_{d\mid a \rad a} G(d)c_d(a) 
=\prod_{p|a} E_{p,G}(a) 
$$
\par
\noindent
(we leave this equation as a kind of exercise for the reader), for which see [{\stampatello Main Lemma}, C2] formul\ae: 
$$
E_{p,G}(a)= \sum_{K=0}^{\infty}G(p^K)c_{p^K}(a)
= \sum_{K=0}^{v_p(a)+1}G(p^K)c_{p^K}(a) 
= \sum_{K=0}^{v_p(a)}p^K\left(G(p^K)-G(p^{K+1})\right), 
$$
\par
\noindent
and (for next equation, see Corollary 3.5 proof) : 
$$
U_{p,G}(a)\defineq G(p^{v_p(a)})-G(p^{v_p(a)+1}) 
\enspace \Rightarrow \enspace 
\UG a \defineq \sum_{d|a}\mu(d)G(da)
=\prod_{p|a} U_{p,G}(a). 
$$
\par
\noindent
Last but not least, we define the other two sums, relative to $G$, for all $a\in\N$: 
$$
\CG a \defineq \sum_{d\mid a} G(d)\mu(d)
=\prod_{p|a}(1-G(p)) 
$$
\par
\noindent
(like for $\EGf$, we leave this equation proof as an exercise [T]), say, the {\it correlating factor} (between $a$ and $G$, see that $\CG a = 0$ iff $a$ has at least one transparent prime factor), and 
$$
\DG a \defineq \sum_{d\mid a} G(d)d, 
$$
\par
\noindent
say, the {\it dilating factor} (for $G$ entails a kind of dilation on $a$ divisors, whenever $|G(p)|>1$ on primes). 

\medskip

We come, of course, to clouds, now. The {\stampatello Ramanujan cloud}, of an arbitrary $F:\N\to\C$ is 
$$
<F>\defineq \{ G:\N\to\C \;|\; \RGf \enspace \hbox{\stampatello is\enspace well-defined\enspace and}\enspace F=\RGf\}, 
$$
\par
\noindent
like defined in [CG], the set of Ramanujan coefficients of this $F$. (See Appendix,A.6: for each $F$, $<F>\neq \emptyset$.) Its most important subset, in this paper, is the {\stampatello multiplicative Ramanujan cloud} of the same $F$ above: 
$$
<F>_M\defineq \{ G\in <F> \;|\; G\enspace \hbox{\stampatello is\enspace multiplicative}\}, 
$$
\par
\noindent
which may be empty; however, for multiplicative $F\neq \0$, $<F>_M\neq \emptyset$, see $\S5$; in which we study the important subset of previous $<F>_M$, that we call the {\stampatello separating multiplicative cloud}, of our $F$, i.e.: 
$$
<F>_{SM}\defineq \{ G\in <F>_M \;|\; \SGf \enspace \hbox{\stampatello is\enspace well-defined\enspace and}\enspace \SG 1\neq 0\}, 
$$ 
\par
\noindent
having in turn the subset which we call the {\stampatello normalized separating multiplicative cloud}
$$
<F>_{NSM}\defineq \{ G\in <F>_M \;|\; \SGf \enspace \hbox{\stampatello is\enspace well-defined\enspace and}\enspace \SG 1=1\}, 
$$ 
\par
\noindent
to which we give a glance in $\S5$. Finally, its notable subset, the {\stampatello Euler separating multiplicative cloud}
$$
<F>_{ESM}\defineq \{ G\in <F>_M \;|\; \SGf \enspace \hbox{\stampatello is\enspace well-defined\enspace and}\enspace \SGf = \1\}, 
$$ 
\par
\noindent
is, for multiplicative $F\neq \0$, a singleton (compare Appendix,A.6), $\{G_F\}$ : Theorem $5.1$ and next definition. 
\par
Last but not least, the construction in Theorem $5.1$ gives the following Ramanujan coefficient. 
\smallskip
\par
\noindent {\bf Definition E.} ({\stampatello Canonical Ramanujan Coefficient of a multiplicative } $F\neq \0$)
$$
G_F(q)\defineq \prod_{p|q}\left(1-\sum_{K=0}^{v_p(q)-1}{{F'(p^K)}\over {p^K}}\right), 
\qquad \forall q\in \N. 
$$

\vfill
\eject

\par				
\noindent{\bf 3. Ramanujan, Lucht and coprime series: from infinite to finite Euler products}
\bigskip 
\par
\noindent
We start quoting a result in our previous paper, [CG], for which we may assume absolute convergence of $\RG a$, for fixed $a$, from absolute convergence of the series of $G(q)\mu(q)$, see Lemma 4 in $\S7.1$. We recognize, in passing, that this is exactly the absolute convergence of $\SG 1$, whence that of every $\SG d$. Then, with this absolute convergence hypothesis, we may write (compare [C2] Proposition): 
$$
\RG a = \prod_{p\in\P}\sum_{K=0}^{v_p(a)}G(p^K)c_{p^K}(a) 
= \prod_{p|a}\sum_{K=0}^{v_p(a)}G(p^K)c_{p^K}(a)\prod_{p\nondivide a}(1-G(p)), 
$$
\par
\noindent
where the Euler product over $p\nondivide a$ is an infinite, but absolutely convergent one. (Compare $\S3$ of [C2].)
\par
This shows clearly that the problem of pointwise convergence of coprime series, without the absolute convergence hypothesis, is a delicate one. Instead, the absolute convergence for $\SGf$ ensures: 
$$
\SG a = \prod_{p\nondivide a}\sum_{K=0}^{v_p(a)}\mu(p^K)c_{p^K}(a) 
= \prod_{p\nondivide a}(1-G(p)), 
$$
\par
\noindent
thus it's a triviality to recover the factorization $\RGf = \EGf \cdot \SGf$ (here $\cdot$ is pointwise product, recall $\EGf$ definition in $\S2$). Also, absolute convergence makes (set $K=tm$ here)
$$
\sum_{K=1}^{\infty}\left|\mu(K)G(dK)\right|\le \sum_{t|d}\sum_{{K=1}\atop {(K,d)=t}}^{\infty}\left|G(dK)\right|\mu^2(K)
\le \sum_{t|d}\mu^2(t)\left|G(dt)\right|\sum_{{m=1}\atop {(m,d)=1}}^{\infty}\left|G(m)\right|\mu^2(m)
<\infty, 
$$
\par
\noindent
ensuring, also, the absolute convergence of $\LG a$, so we may pass to the limit $x\to\infty$ inside: 
$$
\sum_{K\le x}\mu(K)G(dK) = \sum_{t|d}\sum_{{K\le x}\atop {(K,d)=t}}G(dK)\mu(K)
= \sum_{t|d}\mu(t)G(dt)\sum_{{m\le x/t}\atop {(m,d)=1}}G(m)\mu(m)
$$
\par
\noindent
getting : (recall $\UGf$ formul\ae, in $\S2$)
$$
\LG d = \prod_{p|d}\left(G(p^{v_p(d)})-G(p^{v_p(d)+1})\right)\prod_{p\nondivide d}\sum_{K=0}^{\infty}\mu(p^K)G(p^K)
= \UG d \cdot \prod_{p\nondivide d}(1-G(p)), 
$$
\par
\noindent
thus it's again a triviality to factor $\LGf = \UGf \cdot \SGf$. 
\medskip
\par
Of course, absolute convergence (of $\SGf$, meaning : of $\SG 1$, entailing that of all $\SG d$) is a very powerful and useful hypothesis; but we wish not to confine the general frame, so to speak, of our series : in the following, we wish to reach the same conclusions (i.e., factorizations) in the less restrictive assumption, on $\SGf$ series, of pointwise convergence. 
\par
Thus, in the following we will assume convergence of $\SG a$ to factor, first $\RG a$ (next Corollary 3.3) and then $\LG a$ (Corollary 3.5). 
\par
Also, the same philosophy, that we explain in the following (soon before Corollary 3.1), allows to give a kind of \lq \lq recursive formul\ae\rq \rq, for $\SGf$ (Corollary 3.2), for $\RGf$ (Corollary 3.4) and for $\LGf$ (Corollary 3.6). These involve only one series at a time, on which convergence hypotheses, \lq \lq ad hoc\rq \rq, are required; and factors them, according to the corresponding factorization of their argument $a$. This strategy is a precious tool, when proving Theorem 1.1 and Theorem 8.1. 
\par
Last but not least, we'll connect Ramanujan and Lucht series : $G$ is a Ramanujan coefficient (see the above) iff it is a {\stampatello Lucht coefficient}, namely (by definition) $\LG a$ is pointwise converging in all $a\in\N$. (Actually, this equivalence was introduced and proved by Professor Lucht in [Lu1], justifying the names, here.) This is proved in Corollary 3.7; that, actually, proves even more: recalling the Definition D $\S2$ of {\it divisor-closed set}, $\LGf$ converges in a divisor-closed set $\corsivoD$ iff $\RGf$ converges in $\corsivoD$. 

\bigskip

\par				
\noindent
When dealing with an infinite sum of the form \enspace $\sumq h(q)$, where $h:\N\to\C$ is a multiplicative function, a basic technique is to try to write such sum as an infinite Euler product. We can always extract the Euler factors one at at time, to reduce the issues related to convergence: at least formally, we have, for fixed $p$, 
$$
\sumq h(q) = \sum_{k=0}^\infty h(p^k) \sumr p h(r).
$$
\par
\noindent
In fact, this factorization is justified as long as the two factors on the RHS (right-hand side) converge. More generally, we have the following result, which may be considered (like all $\S3$ Corollaries) a \lq \lq consequence\rq \rq, of the arguments proving {\stampatello main lemma} in [C2]: mainly, the {\stampatello vertical limit} \enspace $c_{p^K}(a)=0$, $\forall K>v_p(a)+1$. 
\bigskip
\par
\noindent {\bf Corollary 3.1.} ({\stampatello prototype of finite Euler products extraction})
\par
\noindent
{\it Let $h:\N\to\C$ be a multiplicative function and let $b\in\N$. 
Suppose that  $\sumr b h(r)$ converges and assume that, for each $p|b$, $\exists m_h(p)\in \N:$ $h(p^k)=0, \forall k>m_h(p)$. Then also $\sumq h(q)$ converges and moreover we have
$$
\sumq h(q) = \sum_{\rad d\mid b} h(d) \cdot \sumr b h(r) 
= \prod_{p\mid b} \left( \sum_{k=0}^\infty h(p^k)\right) \cdot \sumr b h(r)
= \prod_{p\mid b} \left( \sum_{k=0}^{m_h(p)} h(p^k)\right) \cdot \sumr b h(r). 
$$
}
\bigskip
\par
\noindent {\it Proof.} Since all $q\in \N$ may be written $q=dr$, where $\rad d|b$ and $(r,b)=1$, entailing $(d,r)=1$, if we set $D\defineq \prod_{p|b}p^{m_h(p)}$, we get 
$$
x>D 
\enspace \Rightarrow \enspace 
\sum_{q\le x}h(q)=\sum_{{d\le x}\atop {\rad d |b}}h(d)\sum_{{r\le x/d}\atop {(r,b)=1}}h(r)
=\sum_{{d\le D}\atop {\rad d |b}}h(d)\sum_{{r\le x/d}\atop {(r,b)=1}}h(r)
=\sum_{\rad d |b}h(d)\sum_{{r\le x/d}\atop {(r,b)=1}}h(r), 
$$
\par
\noindent
so, passing to the limit $x\to \infty$, we obtain first the convergence of LHS (left-hand side) and then, from 
$$
\sum_{\rad d\mid b} h(d)=\prod_{p\mid b} \left( \sum_{k=0}^{m_h(p)} h(p^k)\right)
=\prod_{p\mid b} \left( \sum_{k=0}^\infty h(p^k)\right), 
$$
\par
\noindent
also the formula in the thesis.\qed

\bigskip

\par
\noindent
The convergence of the Euler factor $\sum_{k=0}^\infty h(p^k)$ is, in fact, automatic if $h(p^k)=0$ for all large $k$. This happens for instance if $h$ is supported on the square-free numbers, such as in the case of coprime series, for which $h=G\mu$ ($\mu$ vertical limit's 1) and we have the following.  (Recall $\CGf$ formul\ae, in $\S2$.) 
\smallskip
\par
\noindent {\bf Corollary 3.2.} ({\stampatello The $\SGf-$recursive formula})
\par
\noindent
{\it Let $G:\N\to\C$ be a multiplicative function. Let $b,c\in\N$ be such that $(b,c)=1$. Then
$$
\SGx b x = \sum_{d\mid c} G(d)\mu(d) \SGx {bc} {x/d},
\enspace \forall x\ge 0; 
$$
\par
\noindent
whence, in case $\SG {bc}$ converges, then also $\SG b$ converges and moreover
$$
\SG b = \left(\sum_{d\mid c} G(d)\mu(d)\right) \cdot \SG {bc} 
= \prod _{p|c} (1-G(p)) \cdot \SG {bc}
= \CG c \cdot \SG {bc}. 
$$
}
\par
\noindent {\it Proof.} Since we may write, for any fixed $b\in \N$, all $r\in \N$ with $(r,b)=1$ as $r=du$, with $\rad d|c$ and $(u,bc)=1$, so that $(d,u)=1$, we get, recalling that $x>d$ makes the partial sums \enspace $\SGx {bc} {x/d}$\enspace {\stampatello vanish}: 
$$
\SGx b x = \sum_{{d\mid c}\atop {d\le x}}G(d)\mu(d)\sum_{{(u,bc)=1}\atop {u\le x/d}}G(u)\mu(u)
= \sum_{{d\mid c}\atop {d\le x}}G(d)\mu(d)\SGx {bc} {x/d}
= \sum_{d\mid c}G(d)\mu(d)\SGx {bc} {x/d}, 
\enspace \forall x\ge 0. 
$$
\par
\noindent
Passing to the limit $x\to \infty$, we get second formula.\qed 
\medskip
A similar observation holds for Ramanujan sums, in fact having a vertical limit (like M\"obius function). Then, we get the following formula relating Ramanujan series and coprime series, compare the {\stampatello Proposition} in [C2]. (Recall $\EGf$ formul\ae, in $\S2$.)

\vfill
\eject

\par				
\noindent {\bf Corollary 3.3.} ({\stampatello The $\corsivoR \corsivoS$ formula})
\par
\noindent
{\it Let $G:\N\to\C$ be a multiplicative function and $a\in\N$. Then
$$
\RGx a x =\sum_{d\mid a \rad a} G(d) c_d(a) \SGx a {x/d}, 
\enspace \forall x\ge 0; 
$$
\par
\noindent
whence, in case $\SG a = \SG {\rad a}$ converges, then also $\RG a$ converges and moreover
$$
\RG a = \left(\sum_{d\mid a \rad a} G(d)c_d(a)\right) \cdot \SG a 
= \corsivoE_G(a) \cdot \SG a. 
$$
}
\par
\noindent {\it Proof.} Similarly, any $a\in \N$ admits the factorization $a=dr$, with $\rad d|\rad a$ and $(r,a)=1$, implying again that $(d,r)=1$, so we obtain, again by previous remark on\enspace $\SGf$\enspace partial sums, \enspace $\forall x\ge 0$, 
$$
\RGx a x = \sum_{{d\mid a\rad a}\atop {d\le x}}G(d)c_d(a)\sum_{{(r,a)=1}\atop {r\le x/d}}G(r)\mu(r)
= \sum_{{d\mid a\rad a}\atop {d\le x}}G(d)c_d(a)\SGx a {x/d}
= \sum_{d\mid a\rad a}G(d)c_d(a)\SGx a {x/d}. 
$$
\par
\noindent
Passing to the limit $x\to \infty$ we get second formula.\qed 
\medskip
In order to get next recursive formula for $\RGf$ we apply, together with Corollary 3.2, twice the $\corsivoR \corsivoS$ formula, namely previous Corollary. (Recall $\DGf$ formul\ae, in $\S2$.) 
\smallskip
\par
\noindent {\bf Corollary 3.4.} ({\stampatello The $\RGf-$recursive formula})
\par
\noindent
{\it Fix a multiplicative $G:\N \rightarrow \C$. Let $a,b,c\in\N$ with $a=bc$ and $(b,c)=1$ be such that $\wg=\infty$ for all $p|c$. Then 
$$
\RGx a x = \sum_{h\mid c} G(h) h \enspace \RGx b {x/h},
\enspace \forall x\ge 0. 
$$
\par
\noindent
Consequently, if \thinspace $\RG b$ converges, then \thinspace $\RG a$ converges, as well, and \enspace $\RG a = \DG c \cdot \RG b$. 
}
\smallskip
\par
\noindent {\it Proof.} First, by Corollary 3.3 we have
$$
\RGx a x = \sum_{d\mid a\rad a } G(d) c_d(a) \SGx a {x/d}.
$$
\par
\noindent
For every $d\mid a\rad a$ we may write $d=qr$ where $q \mid b\rad b$ and $r\mid c\rad c$, so 
$$
\RGx a x = \sum_{q\mid b\rad b } G(q) c_q(a)  \sum_{r\mid c\rad c } G(r) c_r(a) \SGx a {x/qr}.
$$
\par
\noindent
We now recall Kluyver's formula (see $\S2$):
$$
c_r(a) = \sum_{h\mid (r,a)} \mu(r/h) h
$$
\par
\noindent
and note that $(r,a)=(r,c)$. We then operate the change of variable $t=r/h$ and make the following change of the order of summation: 
$$
\sum_{r\mid c\rad c} G(r) \sum_{h\mid (r,c)} \mu(r/h) h \SGx a {x/qr} = \sum_{h\mid c} h \sum_{t\mid \rad c} \mu(t)G(th)\SGx a {x/qht}.
$$
\par
\noindent
Since $\wg=\infty$ for all $p|c$, we have that $G$ is completely multiplicative when restricted on the divisors of $c\rad c$. In other words, we have $G(ht)=G(h)G(t)$. Thus, we arrive at the following formula:
$$
\RGx a x = \sum_{q\mid b\rad b } G(q) c_q(a) \sum_{h\mid c} G(h) h \sum_{t\mid \rad c} G(t)\mu(t) \SGx a {x/qht}.
$$
\par				
\noindent
The inner sum is just (recalling: in coprime series, argument's factors, here $c$, can be assumed square-free : $c=\rad c$, here) 
$$
\SGx b {x/qh} = \sum_{t\mid \rad c} G(t)\mu(t) \SGx a {x/qht}, 
$$
\par
\noindent
by Corollary 3.2. Finally, by Corollary 3.3 again, we recognize that
$$
\RGx b {x/h} = \sum_{q\mid b\rad b} G(q) c_q(a) \SGx b {x/qh},
$$
\par
\noindent
and the sought formula is proved.\qed

\bigskip

Our next result is for Lucht series. (Recall $\UGf$ definition, in $\S2$. We prove in next $(\ast)$ its formul\ae.) 
\smallskip
\par
\noindent {\bf Corollary 3.5.} ({\stampatello The $\corsivoL \corsivoS$ formula})
\par
\noindent
{\it Let $G:\N\to\C$ be a multiplicative function. Fix any $d\in\N$. Then 
$$
\LGx d x = \sum_{\ell|d}\mu(\ell)G(\ell d)\SGx d {x/\ell},
\enspace \forall x\ge 0. 
$$
\par
\noindent
Suppose, furthermore, that \enspace $\SG d=\SG {\rad d}$ converges. Then also $\LG d$ converges and moreover
$$
\LG d = \UG d  \cdot \SG d. 
$$
}
\par
\noindent {\it Proof.} Since we may gather, depending on greatest common divisor with any fixed $d\in \N$, 
$$
\LGx d x = \sum_{\ell|d}\sum_{{(K,d)=\ell}\atop {K\le x}}\mu(K)G(dK)
= \sum_{\ell|d}\sum_{{(m,d/\ell)=1}\atop {m\le x/\ell}}\mu(\ell m)G(\ell dm)
= \sum_{\ell|d}\mu(\ell)G(\ell d)\sum_{{(m,d)=1}\atop {m\le x/\ell}}\mu(m)G(m) 
$$
\par
\noindent
we get first formula for all $x\ge 0$; passing to the limit $x\to \infty$ we get second formula. 
\par
We profit here to show that, for multiplicative $G$, we have the following product formula: 
$$
\UG d = \sum_{\ell|d}\mu(\ell)G(\ell d) 
= \prod_{p|d}\left(G(p^{v_p(d)})-G(p^{v_p(d)+1})\right).
\leqno{(\ast)}
$$ 
\par
\noindent
In order to prove this, we first prove that $G$ multiplicative implies $\UGf$ multiplicative. This follows from the definition: 
$$
(a,b)=1
\enspace \Rightarrow \enspace 
\sum_{\ell|ab}\mu(\ell)G(\ell ab)
=\sum_{\ell_a|a}\sum_{\ell_b|b}\mu(\ell_a)\mu(\ell_b)G(\ell_a a)G(\ell_b b)
=\sum_{\ell|a}\mu(\ell)G(\ell a)\cdot \sum_{\ell|b}\mu(\ell)G(\ell b). 
$$
\par
\noindent
Then, the property \lq \lq $\UGf$ multiplicative\rq \rq, just proved, gives the following factorization, that settles $(\ast)$ : 
$$
\sum_{t|d}\mu(t)G(td)=\prod_{p|d}\sum_{\ell|p^{v_p(d)}}\mu(\ell)G(\ell p^{v_p(d)})
=\prod_{p|d}\left(G(p^{v_p(d)}-G(p^{v_p(d)+1}\right), 
\enspace \forall d\in \N.
$$ 
\qed

\bigskip

Our next formula is the recursion for Lucht series. 
\smallskip
\par
\noindent {\bf Corollary 3.6.} ({\stampatello The $\LGf-$recursive formula})
\par
\noindent
{\it Let $G:\N\to\C$ be a multiplicative function and assume $c\in \N$ is such that $p|c$ implies $w_{p,G}=\infty$. Then
$$
\LGx {bc} x = G(c) \LGx b x, 
\enspace \forall x\ge 0. 
$$
\par
\noindent
Suppose, furthermore, that \enspace $\LG b$ converges. Then also $\LG {bc}$ converges and 
$$
\LG {bc} = G(c) \cdot \LG b. 
$$
}

\vfill
\eject

\par				
\noindent {\it Proof.} The assumption on $c$ and $G$ gives $G(cd)=G(c)G(d)$, for all $d\in \N$, whence 
$$
\LGx {bc} x = \sum_{K\le x}\mu(K)G(bcK)
= G(c) \sum_{K\le x}\mu(K)G(bK)
= G(c) \LGx b x, 
\enspace \forall x\ge 0. 
$$
\par
\noindent
When $\LG b$ converges, this entails (passing to $\limx$) convergence of $\LG {bc}$ and 
$$
\LG {bc} = G(c) \cdot \LG b. 
$$
\qed
\medskip
\par
Notice : the property, that we have used many times in previous recursive formul\ae, namely the existence of a factor, say $c\in \N$, having all prime-divisors $p$ such that $w_{p,G}=\infty$, when considered in the case of bad (see $\S2$ definitions) primes $p|c$, is equivalent in saying that $c$ is made of hyperbad (see $\S2$ again) prime factors; this, in turn, see the definition of Ramanujan factorization of a natural number in $\S8.1$, amounts to saying that $h(c)=c$. In fact, the recursive formul\ae, as we will see, are vital when proving Theorem 8.1. 
\medskip
Now, we may give our last Corollary; that is, also, the last possible connection in between $\SGf,\RGf,\LGf$ (including the recursions) : the $\corsivoL \corsivoR$ link. A great difference with ALL previous formul\ae, here, is the generality: $G$ is any arithmetic function. We are grateful to Professor Lucht. 
\smallskip
\par
\noindent {\bf Corollary 3.7.} ({\stampatello The $\corsivoL \corsivoR$ formula}) {\it Let } $G$ {\it be} {\stampatello any arithmetic function.} {\it Then} 
$$
\LG d \enspace \hbox{\stampatello converges}\enspace  \forall d|a 
\enspace \Rightarrow \enspace 
\RG a \enspace \hbox{\stampatello converges}
$$
\par
\noindent
{\it and 
$$
\RG d \enspace \hbox{\stampatello converges}\enspace \forall d|a 
\enspace \Rightarrow \enspace 
\LG a \enspace \hbox{\stampatello converges}, 
$$
\par
\noindent
whence } 
$$
G \enspace \hbox{\stampatello is\enspace a\enspace Ramanujan\enspace coefficient}
\enspace \Longleftrightarrow \enspace 
G \enspace \hbox{\stampatello is\enspace a\enspace Lucht\enspace coefficient}. 
$$
\par
\noindent
{\it In either case, 
$$
\RG a = \sum_{d|a}d\;\LG d, 
\enspace \forall a\in \N
\quad \hbox{\stampatello and} \quad
\LG d = {1\over d} \sum_{t|d}\RG t \mu(d/t),
\enspace \forall d\in \N. 
$$
}
\par
\noindent {\it Proof.} Clearly, it suffices, here, to prove only the following two partial sums formul\ae. 
\par
\noindent
First one follows from writing, from Kluyver's formula (see $\S2$): 
$$
\RGx a x = \sum_{d|a}d\sum_{{q\le x}\atop {q\equiv 0\bmod d}}\mu(q/d)G(q)
= \sum_{d|a}d\sum_{K\le x/d}\mu(K)G(dK)
= \sum_{d|a}d \; \LGx d {x/d},
\quad \forall x\ge 0; 
$$
\par
\noindent
whence, by M\"obius inversion, see [T], also second formula follows: 
$$
\LGx d {x/d} = {1\over d}\sum_{t|d}\RGx t x \mu(d/t), 
\quad \forall x\ge 0. 
$$
\qed

\bigskip

\par
This Corollary regards the wonderful link between Ramanujan series and Lucht series; actually, it is a kind of writing explicitly an idea going back to Lucht [Lu1] (and compare [Lu2] for a better explanation), namely of using the Eratosthenes transform (see $\S2$ \& $\S4$) to connect these two. It's because of this, that we called Lucht series the series introduced in this paper. Actually, see that it is exactly the last parts of previous Corollary, that were proved by Lucht; however, his elementary argument proves even the formul\ae, \thinspace on divisors, of this result. By these formul\ae: $\RGf$ is {\stampatello well-defined on} a divisor-closed $\corsivoD$ (Definition D,$\S2$) {\stampatello if and only if} $\LGf$ is {\stampatello well-defined on} $\corsivoD$. 
\par
We'll apply this last Corollary in our next section for completely multiplicative clouds. 

\vfill
\eject

\par				
\noindent{\bf 4. Eratosthenes transforms and completely multiplicative clouds}
\bigskip
\par
\noindent
In the following, the Completely Multiplicative Ramanujan cloud of a fixed $F:\N\to\C$ is 
$$
<F>_{CM}\defineq \left\{G\in <F> \left.\right| G\; \hbox{\rm is} \enspace \hbox{\stampatello completely \enspace multiplicative}\right\}. 
$$
\par
Recall the Eratosthenes Transform $f'$ of any $f:\N\to\C$, in Definition C of $\S2$, and that $f'$ is multiplicative iff this $f$ is multiplicative. 
\par
From Corollary 3.7 above, for all $G:\N\to\C$ we may calculate $\corsivoR_G'$, the Eratosthenes Transform of \enspace $\RGf$: 
$$
\RG a =\sum_{d|a}\ETRG d, 
\enspace \forall a\in\N
\quad \Rightarrow \quad
\ETRG d = d\;\LG d, 
\enspace \forall d\in\N,
$$
\par
\noindent
showing a r\^ole of Lucht series, as a kind of (normalized) Eratosthenes transform of Ramanujan series. 
\medskip
\par
The Lucht series of a completely multiplicative coefficient $G$ is very easy. 
\smallskip
\par
\noindent {\bf Lemma 4.1.} {\it Let \enspace $G:\N\to\C$ be completely multiplicative, for which } $\LG 1$\enspace {\it converges. Then} 
$$
\LG a = G(a) \LG 1,
\quad \forall a\in \N. 
$$
\par
\noindent {\it Proof.} It follows from $\LGf$ definition.\qed   
\medskip
\par
We describe $<F>_{CM}$ now. We write \lq \lq $C.M.$\rq \rq, for \lq \lq is completely multiplicative\rq \rq. 
\smallskip
\par
\noindent {\bf Theorem 4.2.} {\it Let \enspace $F:\N\to\C$ be any arithmetic function, not the null-function (i.e., $F\neq \0$). Then 
$$
<F>_{CM}\neq \emptyset
\quad \Longleftrightarrow \quad 
\sum_{q=1}^{\infty}{{F'(q)}\over q}\mu(q) \; \hbox{\it converges\enspace to} \enspace F^2(1)\neq 0 
\enspace \hbox{\stampatello and} \enspace 
F'/F(1) \enspace \hbox{\it C.M.} 
$$
\par
\noindent
When this is the case, $<F>_{CM}$ is the singleton \enspace 
${\displaystyle 
\left\{ G:\N\to\C\; |\; G(q)={{F'(q)}\over {qF(1)}},\, \forall q\in \N\right\}. 
}$
}
\smallskip
\par
\noindent {\it Proof.} First of all, hereafter $\RG 1 = F(1)\neq 0$; otherwise, for all $G$ C.M., next Lemma gives $\RGf = F = \0$. 
\par
We prove, first, the \thinspace $(\Rightarrow)$. 
From Corollary 3.7, $G\in<F>$ $\Rightarrow$ $F'(d)=\corsivoR_G'(d)=d\,\LG d$ and then from Lemma 4.1, $G$ C.M. $\Rightarrow$ $\LG d = G(d)\LG 1 = G(d)F(1)$, whence $F'(d)=dG(d)F(1)$, $\forall d\in\N$.\hfill QED 
\par
While, \enspace $(\Leftarrow)$ \enspace follows from choosing $G(d)={{F'(d)}\over {dF(1)}}$, that's C.M. from $F'/F(1)$ C.M., so that following Lemma with $F'$ convergence hypothesis give the convergence of 
$$
\RG a = \DG a \cdot \sum_{q=1}G(q)\mu(q) 
 = \sum_{d|a}{{F'(d)}\over {F(1)}} \cdot \sum_{q=1}{{F'(q)}\over {F(1)q}}\mu(q)
 = \sum_{d|a}F'(d)
 = F(a),
\enspace \forall a\in\N, 
$$ 
\par
\noindent
whence : $G\in<F>_{CM}$.\qed 
\medskip
\par
\noindent {\bf Remark 1.} {\it See that $F'/F(1)=(F/F(1))'$, by linearity of Eratosthenes transform, so, if $F'/F(1)$ is C.M., in particular it's multiplicative and then: $F'/F(1)=(F/F(1))'$ is multiplicative, implying (Eratosthenes Transform and its inverse preserve multiplicativity) $F/F(1)$ multiplicative, whence $F$ is quasi-multiplicative. An arithmetic function $F$ is quasi-multiplicative, in fact, by definition, iff there's a constant $c\neq 0$ (to avoid trivialities), with $F/c$ multiplicative; and, in turn (avoiding trivial cases), this implies $c=F(1)$. For this definition and properties of quasi-multiplicative functions, see $[{\rm Lah}]$.}\hfill $\diamond$ 
\medskip
\par
This Theorem allows us to describe completely, once given any $F:\N\to\C$ (which is not the null-function, otherwise see Theorem 9.1 and compare [CG]), its C.M. Ramanujan cloud $<F>_{CM}$; in particular, depending on $F$, it may only be empty or a singleton, a kind of \lq \lq uniqueness\rq \rq. 
\par
Theorem 4.2 follows the general philosophy of \lq \lq uniqueness results\rq \rq, for Ramanujan expansions. As explained in [C3], if, for a fixed $F:\N\to\C$, we want a {\bf unique} Ramanujan expansion, from the features of this $F$, this {\bf is} {\bf impossible}: we can always add to a fixed Ramanujan coefficient $G:\N\to\C$ of our $F$, say, $G_0\in<\0>$, i.e. any coefficient of null-function $\0$. Instead, we have to impose conditions on $G$, not on $F$ (compare [C3, $\S3$], Theorems 1,2,3). In Theorem 4.2, the condition is: $G$ is completely multiplicative. 

\bigskip

\par
Last but not least, the Ramanujan series with completely multiplicative coefficient $G$ is very easy, too. 
\smallskip
\par
\noindent {\bf Lemma 4.3.} {\it Let \enspace $G:\N\to\C$ be completely multiplicative, for which } $\RG 1$\enspace {\it converges. Then} 
$$
\RG a = \DG a \RG 1, 
\quad \forall a\in \N. 
$$
\par
\noindent {\it Proof.} It follows from the $\RGf-$recursive formula, i.e. Corollary 3.4, taking $b=1$ and $c=a$.\qed   

\vfill
\eject

\par				
\noindent{\bf 5. Canonical Ramanujan coefficients and the separating multiplicative clouds}
\bigskip
\par
\noindent
We start, in order to study the \lq \lq {\stampatello separating multiplicative clouds}\rq \rq, recalling their definition (in $\S2$): 
$$
<F>_{SM}\defineq \{G\in <F>_M\, :\, \SGf\, \hbox{\stampatello well-defined\enspace and}\enspace \SG 1\neq 0\}, 
$$
\par
\noindent
which we don't discuss, for the moment. In fact, even the reason why we chose the name \lq \lq separating\rq \rq \thinspace will be explained later. 
\smallskip
\par
Our strategy will be to start with \lq \lq simpler\rq \rq, say, subsets of $<F>_{SM}$; actually, these are much smaller than it (but this will be clear a posteriori) and, in fact, we'll start with the smallest, following; then, we'll take a glance at an intermediate one; and, finally, at last in this section, we'll explain (the name origin and) the properties of $<F>_{SM}$. Luckily enough, we can introduce the {\it canonical Ramanujan coefficients}, with the smallest, $<F>_{ESM}$, following. 
\smallskip
\par
We start, in order to arrive to the study of separating multiplicative clouds, with the following subset, for a fixed $F:\N\to\C$
$$
<F>_{ESM}\defineq \{G\in <F>_M\, :\, \SGf=\1\} 
$$
\par
\noindent
(notice that $F(1)=\SG 1=1$ excludes $F=\0$, namely, $<\0>_{ESM}=\emptyset$), which is the {\stampatello Euler separating multiplicative cloud} of this $F$ (compare the following, for the reason why the name \lq \lq separating\rq \rq). For a fixed $F:\N\to\C$, assumes that the coprime series is the $\1$ function (constantly $1$ on $\N$, $\S2$), so that, from the $\corsivoR \corsivoS$ formula (Corollary 3.3), we get $\RGf=\EGf$, whence name \lq \lq Euler\rq \rq: the Ramanujan series, here, is the same Euler-Ramanujan factor. In other words, $F$ has a finite Ramanujan expansion, with explicit finite Euler product and is multiplicative, because $\EGf=\RGf=F$ is multiplicative ! Needless to say, for $F$ not multiplicative, $<F>_{ESM}=\emptyset$. 
\par
Furthermore, $G\in <F>_{ESM}$ implies, in particular, that $G$ vanishes on primes (which is easily proved by $iii)$ in Lemma 5.3, holding in $<F>_{SM}$); more in general, when $G$ is not multiplicative, we say that 
$$
G:\N\to\C
\enspace \hbox{\it is} \enspace \hbox{\stampatello square-free\enspace vanishing} 
\quad \definiz \quad
G(s)=0, \enspace \forall s>1, \mu^2(s)=1
$$
\par
\noindent
which can be stated, in fact, for multiplicative $G$ as : $G(p)=0$ for all primes $p$. 
\par
Then we are, actually, building one possible multiplicative Ramanujan coefficient, of a multiplicative, fixed $F:\N\to\C$. We know it has to vanish on primes; but we still don't know how it's made on prime-powers: notice that this information is all built-in the Euler-Ramanujan factor ! An explicit  construction, for the {\stampatello canonical Ramanujan coefficient} $G_F$ ($\S2$) of any multiplicative $F\neq \0$, is given in next result proof. 
\medskip
\par 
Thus any Euler separating multiplicative cloud of a multiplicative $F\neq \0$ is a singleton. 
\smallskip
\par
\noindent {\bf Theorem 5.1.} {\it Let\enspace $F\neq \0$ \enspace be {\stampatello multiplicative}. Then } $<F>_{ESM}=\{G_F\}$, {\it whence } $<F>_M\enspace \neq \enspace \emptyset$. 
\smallskip
\par
\noindent {\it Proof.} Once fixed $F:\N \rightarrow \C$, multiplicative, since we exclude $F=\0$($\defineq$the null-function, see [CG]), hereafter $a=1$ has empty product on $p|a$, with value $1$ in the following : 
$$
F(a)=\prod_{p|a}F(p^{v_p(a)})
\quad \forall a\in \N; 
$$
\par
\noindent
and we'll build a multiplicative $G:\N \rightarrow \C$, say, $G=G_F$ ($G$ depends on $F$), then proving $<F>_{ESM}=\{G_F\}$.
\par
\noindent
We start from the hypothesis that $G$ is multiplicative with $G(1)=1$ (to avoid $G=\0$) and square-free vanishing. Then, (see the above) 
$$
\RG a=\sum_{q=1}^{\infty}G(q)c_q(a)
=\EG a
=\prod_{p|a}\sum_{K=0}^{v_p(a)}p^K\left(G(p^K)-G(p^{K+1})\right). 
$$
\par
\noindent
Up to now we have only imposed $G_F$ multiplicative and $G_F(1)\defineq 1$, $G_F(p)\defineq 0$ on all primes $p$. (So we know $G_F\in <F>_{ESM}$.) 
\par				
\noindent
Thus, one natural step is to check if $\RG a = F(a)$, for all natural $a$ (both sides are defined globally); so, we make, once fixed $p|a$, 
$$
F(p^{v_p(a)})=\sum_{K=0}^{v_p(a)}p^K\left(G(p^K)-G(p^{K+1})\right)
$$
\par
\noindent
an assumption that, made $\forall a\in \N$, with $p|a$, is equivalent to the following: 
$$
F(p^v)=\sum_{K=0}^{v}p^K\left(G(p^K)-G(p^{K+1})\right),
\quad \forall v\in \N. 
\leqno{(\ast)_{G_F}}
$$
\par
\noindent
Let's see our first definition, setting $v=1$ in $(\ast)_{G_F}$ and using $G(1)=1, G(p)=0$ : 
$$
G_F(p^2)\defineq {{1-F(p)}\over {p}}, 
\quad \forall p\in \P. 
$$
\par
\noindent
More in general, by induction on $v\in \N$, equation $(\ast)_{G_F}$ above allows us to build $G_F$ on some power, once defined on lower powers ($K=1$ has square-free vanishing since $K=0$ forces $G_F(1)$ to be $1$). 
\par
\noindent
Calculating $F(p^v)-F(p^{v-1})$ from $(\ast)_{G_F}$, we get $G(p^v)-G(p^{v+1})$, whence, in fact, by induction and also abbreviating $F(p^v)-F(p^{v-1})=F'(p^v)$, where $F'$ is the Eratosthenes transform [Wi] of our $F$ (see $\S2$, $\S4$), we define (here for $v=1$ the empty sum over $1\le K\le v-1$ gives, in fact, $G_F(p)=0$)
$$
G_F(p^v)\defineq -\sum_{K=1}^{v-1}{{F'(p^K)}\over {p^K}}
=1-\sum_{K=0}^{v-1}{{F'(p^K)}\over {p^K}}, 
\qquad \forall v\in \N
$$
\par
\noindent
and, recalling: $f$ multiplicative $\Rightarrow $ $\prod_{p|n}(1-f(p))=\sum_{t|n}\mu(t)f(t)$, with $f(t)=\prod_{p|t}f(p)$, we give, compare $\S2$ definition, the formula (entailing $G_F(1)=1$, from: empty products are $1$ and $G_F(p)=0$, from: $v_p(p)=1$)
$$
G_F(q)=\sum_{t|q}\mu(t)\prod_{p|t}\sum_{K=0}^{v_p(q)-1}{{F'(p^K)}\over {p^K}}, 
\qquad \forall q\in \N. 
$$
\par
\noindent
Thus
$$
G=G_F \thinspace \hbox{\rm inside} \thinspace (\ast)_{G_F}
\thinspace \Rightarrow \thinspace 
F(a)=\prod_{p|a}F(p^{v_p(a)})=\prod_{p|a}\sum_{K=0}^{v_p(a)}p^K\left(G(p^K)-G(p^{K+1})\right)
=\EG a=\RG a, \forall a\in \N, 
$$
\par
\noindent
having on $a=1$ the convention of empty products being $1$. This procedure proves $<F>_{ESM}=\{G_F\}$.\qed 
\medskip
\par
Hereafter we need a name for this $G_F$, since it's the {\bf unique Ramanujan coefficient} inside $<F>_{ESM}$: we call it (see $\S2$) the {\stampatello canonical Ramanujan coefficient} of $F$; enlarging this set, uniqueness is lost. 
\par
We relax the condition $\SGf=\1$ with a more general one, namely $\SG 1=1$, getting the following {\stampatello normalized separating multiplicative clouds}: 
$$
<F>_{NSM}\defineq \{G\in <F>_M\, :\, \SGf\, \hbox{\stampatello well-defined\enspace and}\enspace \SG 1=1\}, 
$$
\par
\noindent
which is, still, a subset of $<F>_{SM}$, better, $<F>_{ESM}\subset <F>_{NSM}\subset <F>_{SM}$ (and $N$ stands for \lq \lq normalized\rq \rq, somehow imagining to divide by $\SG 1\neq 0$). A simple, but discouraging, remark is that, now, a priori, we may loose $F$ multiplicativity ! However, from Corollary 1.2 and the $G$ multiplicativity, we always have a semi-multiplicative $F$ : then, noticing $F(1)=\SG 1=1$, we recover that $F$ is multiplicative (since has threshold $a_F=1$ and $c=F(a_F)=1$) ! Even better, \thinspace $<F>_{NSM}\neq \emptyset$\enspace $\Longleftrightarrow$\enspace $F$ is multiplicative. 
\par
Actually, a kind of \lq \lq raison d'\^{e}tre\rq \rq, for this subset, is just to keep the property: $F$ multiplicative. In fact, we'll not give, for the time being, any further property of $<F>_{NSM}$. Time to start the study of $<F>_{SM}$ now, giving a kind of \lq \lq coming soon\rq \rq: this time (in order to avoid $<F>_{SM}=\emptyset$) $F$ is quasi-multiplicative (more general than \lq \lq $F$ is multiplicative\rq \rq, compare the definition in Remark 1, soon after Theorem 4.2). 
\medskip
\par
One good idea is to start from the reason why: \lq \lq separating\rq \rq. Actually, we thought to mean: this $<F>_{SM}$ somehow separates, for each $\SG a$, the primes dividing $a$ from all the others and it will be clear (compare Lemma 5.3) from the finite product (and infinite product, too) formul\ae, coming next. However, other properties, in Lemma 5.3, may justify the name as well. 

\vfill
\eject

\par				
\noindent
Before we give the important Lemma 5.3, encoding the main features of separating multiplicative coefficients, we wish to underline (in next Lemma) the importance of absolute convergence, for coprime series, in order to simplify all the formul\ae, giving infinite products. 
\medskip
\par
Recall $G$ is normal iff, by definition, it's multiplicative and has no transparent primes, see $\S2$ and [CG]. 
In the following, we'll use the adjective \lq \lq separating multiplicative\rq \rq, not only for clouds, but also for arithmetic functions (typically, $G$). 
\medskip
\par
The normal arithmetic functions $G$ with an absolute convergence condition are separating. 
\smallskip
\par
\noindent {\bf Lemma 5.2.} ({\stampatello Normal $G$ with absolutely convergent $\SG 1$ are separating multiplicative})
\par
\noindent
{\it Let } $G:\N\to\C$ {\it be normal, with } $\SG 1$ {\it absolutely convergent. Then 
$$
\SGf\enspace \hbox{\stampatello converges\enspace absolutely\enspace and is never-vanishing}; 
$$
\par
\noindent
in formul\ae,
$$
\SG a = \prod_{p\nondivide a}(1-G(p))\neq 0\enspace \hbox{\stampatello converges\enspace absolutely}
\quad \forall a\in \N; 
$$
\par
\noindent
in particular, $G$ is {\stampatello separating multiplicative}.
}
\smallskip
\par
\noindent {\it Proof.} The absolute convergence of \enspace $\SG 1$, by positivity (recall \thinspace $\SGf$ \thinspace definition, $\S1$), implies that of any \enspace $\SG a$, at once. 
\par
\noindent
The passage from $\SG a$ absolute convergence to the absolute convergence of its Euler product above, then, is standard (details in [C2, $\S3$]).  
\par
\noindent
A classical result (dating back to \lq \lq old times\rq \rq), then, ensures that an absolutely convergent infinite product which has no vanishing factor (as $G$ is  normal here) doesn't vanish (details, again, in [C2, \S3], quoting Ahlfors Book). 
\qed 
\medskip
\par
We come to the main features of coprime series, when $G$ is separating multiplicative. 
\smallskip
\par
\noindent {\bf Lemma 5.3.} ({\stampatello Coprime series with separating multiplicative } $G$)
\par
\noindent
{\it Let } $G:\N\to\C$ {\it be} {\stampatello multiplicative, with $\SGf$ well-defined and $\SG 1\neq 0$.} {\it Then} 
\medskip
\par
\noindent
$i)$ {\it $G$ is} {\stampatello normal;}
$$
\SG a = \prod_{p|a}(1-G(p))^{-1}\, \SG 1
= \sum_{d|a}{{\mu(d)G(d)}\over {\prod_{p|d}(G(p)-1)}}\SG 1,
\quad
\forall a\in \N;
\leqno{ii)}
$$
\par
\noindent
$iii)$ {\it in particular, } $G(p)=1-\SG 1/\SG p$, $\forall p\in\P$, {\it so } $\SGf$ {\it on $\{1\}\cup \P$ determines } $G$ {\it on } $\P$; 
\medskip
\par
\noindent
$iv)$ $\SG 1$ {\stampatello converges absolutely\enspace iff \enspace } $\sum_{p\in\P}\left|1-\SG 1/\SG p\right|<\infty$;
\medskip
\par
\noindent
$v)$ $\SGf$ {\stampatello is never-vanishing}. 
\smallskip
\par
\noindent {\it Proof.} The $\SGf-$recursive formula (Corollary 3.2) gives 
$$
\SG a \prod_{p|a}(1-G(p))=\SG 1,
\quad 
\forall a\in \N, 
$$
\par
\noindent
whence \enspace $\SG 1 \neq 0$ \enspace proves $i)$, then $iii)$ and the first equation in $ii)$. Since (compare Theorem 5.1 proof) 
$$
\prod_{p|a}(1-G(p))^{-1}=\prod_{p|a}{1\over {1-G(p)}}
=\prod_{p|a}\left(1-{{G(p)}\over {G(p)-1}}\right)
=\sum_{d|a}\mu(d)\prod_{p|d}{{G(p)}\over {G(p)-1}},
\quad 
\forall a\in \N, 
$$
\par
\noindent
we get, as $\mu(d)\neq 0$ implies $\prod_{p|d}G(p)=G(d)$, the second equation in $ii)$. 
\par
The $\SG 1$ absolute convergence, thanks to Lemma 4 in $\S7.1$ of [CG], is equivalent to $\sum_{p\in\P}|G(p)|<\infty$; then, we conclude $iv)$ proof, by $iii)$, proved above. 
\par
Finally, $v)$ follows immediately, from $ii)$ and the hypothesis $\SG 1\neq 0$.
\qed
\medskip
\par				
\noindent {\bf Remark 2.} {\it Hence, abbreviating {\stampatello w-d for well-defined} (recall, pointwise convergence in all\enspace $\N$)
$$
\SGf\enspace \hbox{\stampatello w-d\enspace and}\enspace \SG 1\neq 0
\quad \Longleftrightarrow \quad
\SGf\enspace \hbox{\stampatello w-d\enspace and\enspace never-vanishing}, 
$$
\par
\noindent
rendering RHS condition an equivalent, possible alternative definition of \lq \lq $G$ separating multiplicative\rq \rq.}\hfill $\diamond$
\medskip
\par
We conclude with a property regarding separating multiplicative clouds i.e. $<F>_{SM}\neq \emptyset$ implies that, in general (since $\SG 1\neq 0$) $F$ is quasi-multiplicative (see Remark 1), while adding the condition $\SG 1=1$, i.e. $<F>_{NSM}\neq \emptyset$, implies that $F$ is multiplicative (see the above). 
\par
The quasi-multiplicativity alone, for the $\RGf$, in case the multiplicative Ramanujan coefficients $G$ have $\SGf$ well-defined (in [La], see Definition 1.2) and $\RG 1 \neq 0$ (equivalent to $\SG 1 \neq 0$), namely for separating multiplicative $G$, was already recognized in Corollary 1.2 of [La] (proved with same ideas, different details). 
\medskip

\vfill
\eject

\par				
\noindent{\bf 6. Converse convergence theorem}
\bigskip 
\par
\noindent
We may wonder if a converse of Corollary 3.2 holds, that is, if the convergence of $\SG b$ for some $b\in\N$ implies the convergence of $\SG {bc}$ for all $c\in\N$ coprime with $b$. This is false in general: we give some counterexamples in Appendix, A.2. However, it is true as long as $c$ has no prime divisors with $1\leq |G(p)|\leq p$. We prove this result in Corollary 6.2 below. The key tool is the following \lq \lq converse convergence theorem\rq \rq. We remark that the results of this  section are the motivation for the introduction of the notion of \lq \lq bad primes\rq \rq \thinspace in Section 2. 
\par
We see the result proving Corollary 6.2; however, it also has an interest of his own. 
\smallskip
\par
\noindent {\bf Theorem 6.1.} {\stampatello (Converse convergence theorem)}
\par
\noindent
{\it Let $\alpha\in\C\backslash \{-1\}$ and $\rho>1$. Assume that the function \enspace $H:\R_{\ge 0}\rightarrow \C$ \enspace satisfies $H(x)=H([x])$, \enspace $\forall x\ge 0$ \enspace and \enspace $K(x) := H(x) + \alpha H(x/\rho)\to \ell \in \C$, as $x\to\infty$. Then, $\limx H(x)$ exists in $\C$ and equals $\ell/(1+\alpha)$, in both of the distinct hypotheses }
$$
|\alpha|<1; 
\leqno{\hbox{\stampatello (contraction)}}
$$
$$
|\alpha|>\rho \quad \hbox{\it and} \quad H(x)\ll x, \enspace \hbox{\it as} \enspace x\to \infty.
\leqno{\hbox{\stampatello (dilation)}}
$$ 
\medskip
\par 
\noindent {\bf Remark 3.} {\it We may relax the hypothesis on $H([x])$, assuming only $H$ bounded on $\R_{\ge 0}$ compact subsets.}\hfill $\diamond$ 
\medskip
\par
\noindent {\it Proof.} First, the \lq \lq contractive case\rq \rq, \thinspace $0<|\alpha|<1$ (since $\alpha=0$ is trivial). 
\par
\noindent
We start proving that $H$ is bounded on $\R_{\ge 0}$. Then, this will imply $\limx H(x)=\ell/(1+\alpha)$. 
\par
\noindent
Since $K(x)$ has a finite limit, as $x\to \infty$, and $H$ (from hypothesis on $[x]$, above) is bounded on compact subsets of $\R_{\ge 0}$, we derive
$$
M:=\sup_{x\ge 0}|H(x) + \alpha H(x/\rho)|<\infty. 
$$
\par
\noindent
The triangle inequality gives from this
$$
|H(x)|=|H(x)+\alpha H(x/\rho)-\alpha H(x/\rho)|
       \le M + |\alpha|\cdot |H(x/\rho)|,
\quad \forall x\ge 0, 
$$
\par
\noindent
whence, iterating,
$$
|H(x)|\le M (1+|\alpha|+\cdots+|\alpha|^{j-1}) + |\alpha|^{j}\cdot |H(x/\rho^{j})|,
\quad \forall x\ge 0, \forall j\in \N. 
$$
\par
\noindent
Fix $n\in \N$ and $0\le x\le n$. Completing to the geometric series of ratio $|\alpha|$ and using $\rho>1$, 
$$
|H(x)|\le {M\over {1-|\alpha|}} + |\alpha|^{j}\cdot \sup_{0\le x\le n}|H(x/\rho^{j})|
       \le {M\over {1-|\alpha|}} + |\alpha|^{j}\cdot \sup_{0\le x\le n}|H(x)|,
\quad \forall j\in \N. 
$$
\par
\noindent
Then, once fixed $\varepsilon>0$, we may choose $j\in \N$, depending on $\varepsilon$ and $n,H$, (esp., $j>(\log{{J(n)}\over {\varepsilon}})/(-\log|\alpha|)$, for $J(n):=\sup_{0\le x\le n}|H(x)|$, of course $J(n)>0$ here) such that 
$$
0\le x\le n 
\enspace \Longrightarrow \enspace 
|H(x)|\le {M\over {1-|\alpha|}} + \varepsilon. 
$$
\par
\noindent
Since this is true $\forall n\in \N$ (and notice that $\varepsilon$ doesn't depend on $n$), $H$ is bounded on all $\R_{\ge 0}$.
\hfill (QED)
\par
\noindent
Once proved this, we apply it now, assuming $S_H:=\sup_{x\ge 0}|H(x)|$ positive to avoid trivialities. 
\par
\noindent
From the hypothesis, for all $n\in \N$: 
$$
(-\alpha)^j(H(x/\rho^j)+\alpha H(x/\rho^{j+1}))\to (-\alpha)^j\ell, 
\qquad \forall j=0,\ldots,n-1; 
$$
\par
\noindent
summing up, 
$$
\limx \left(H(x)-\alpha^n H(x/\rho^n)\right)={{1-(-\alpha)^n}\over {1+\alpha}}\ell
={{\ell}\over {1+\alpha}}-{{(-\alpha)^n}\over {1+\alpha}}\ell, 
$$
\par				
\noindent
i.e., from the definition of limit, 
$$ 
\forall \varepsilon>0 \  \exists x_{\varepsilon}\ge 0
\enspace : \enspace 
\left|H(x)-\alpha^n H(x/\rho^n)-\left({{\ell}\over {1+\alpha}}-{{(-\alpha)^n}\over {1+\alpha}}\ell\right)\right|<{{\varepsilon}\over 2},
\quad \forall x>x_{\varepsilon} 
$$
\par
\noindent
and we may choose $n\in \N$ so large that (esp., $n>(\log{{2C}\over {\varepsilon}})/(-\log|\alpha|)$, writing $C:=|S_H|+{{|\ell|}\over {|1+\alpha|}}$)
$$ 
\left|-\alpha^n H(x/\rho^n)+{{(-\alpha)^n}\over {1+\alpha}}\ell\right|
\le |\alpha|^n\cdot \left(|S_H|+{{|\ell|}\over {|1+\alpha|}}\right)
<{{\varepsilon}\over 2},
\quad \forall x\ge 0, 
$$
\par
\noindent
whence the triangle inequality gives immediately $H(x)\to {{\ell}\over {1+\alpha}}$, from the definition of limit.
\hfill (QED)
\par
\noindent
We prove, as above, first that $H$ is bounded on $\R_{\ge 0}$. For this, we mimick previous argument, giving less details, of course. Now, the main change will be using ratio $1/\alpha$, instead of $\alpha$. Also, \hfil $1/|\alpha|$ \hfil substitutes: \hfil $|\alpha|$. 
\par
\noindent
As above, $K(x)\to \ell$ finite and $H$ is bounded on compacta of $\R_{\ge 0}$, so we get ($L=M/|\alpha|$, compare above): 
$$
L:=\sup_{x\ge 0}|H(\rho x)/\alpha + H(x)|<\infty. 
$$
\par
\noindent
The triangle inequality, iteration and completing to the geometric series give, once fixed $n\in\N$,  
$$
|H(x)|\le L + {{|H(\rho x)|}\over {|\alpha|}}, 
\enspace \forall x\ge 0 
\enspace \Rightarrow \enspace 
|H(x)|\le {L\over {1-|\alpha|^{-1}}} + {1\over {|\alpha|^{j}}}\cdot \sup_{0\le x\le n}|H(\rho^{j} x)|, 
\enspace \forall j\in \N. 
$$
\par
\noindent
The hypotheses on $H$ give that $\exists C,D>0$ (absolute constants): $|H(x)|\le Cx+D$, $\forall x\ge 0$, entailing 
$$
|H(x)|\le {L\over {1-|\alpha|^{-1}}} + {{C\rho^{j} n+D}\over {|\alpha|^{j}}}
       \le {L\over {1-|\alpha|^{-1}}} + Cn\left({{\rho}\over {|\alpha|}}\right)^{j}+{D\over {|\alpha|^j}},
\quad \forall j\in \N. 
$$
\par
\noindent
Hence, once fixed $\varepsilon>0$, we may choose $j\in \N$, depending on $\varepsilon$ and $n,C,D$, so large that 
$$
0\le x\le n 
\enspace \Longrightarrow \enspace 
|H(x)|\le {L\over {1-|\alpha|^{-1}}} + \varepsilon. 
$$
\par
\noindent
Since this is true $\forall n\in \N$ (and notice that $\varepsilon$ doesn't depend on $n$), $H$ is bounded on all $\R_{\ge 0}$.\hfill (QED)
\par
\noindent
Assume again previous $S_H>0$. From the hypothesis on $K(x)\to \ell$, for all $n\in \N$, summing up on $0\le j<n$,  
$$
\left(-{1\over {\alpha}}\right)^j\left(H(\rho^j x)+{1\over {\alpha}}H(\rho^{j+1} x)\right)\to {{\ell}\over {\alpha}}\left(-{1\over {\alpha}}\right)^j 
\enspace \Rightarrow \enspace
H(x)-\left(-{1\over {\alpha}}\right)^n H(\rho^n x)\to {{\ell}\over {1+\alpha}}-{{\ell}\over {1+\alpha}}\left(-{1\over {\alpha}}\right)^n 
$$
\par
\noindent
i.e., from the definition of limit, 
$$ 
\forall \varepsilon>0 \  \exists x_{\varepsilon}\ge 0
\enspace : \enspace 
\left|H(x)-{{\ell}\over {1+\alpha}}+\left({{\ell}\over {1+\alpha}}-H(\rho^n x)\right)\left(-{1\over{\alpha}}\right)^n\right|<{{\varepsilon}\over 2},
\quad \forall x>x_{\varepsilon} 
$$
\par
\noindent
and we may choose $n\in \N$ so large that 
$$ 
\left(\left|{{\ell}\over {1+\alpha}}\right|+\left|H(\rho^n x)\right|\right){1\over{|\alpha|^n}}
\le \left(\left|{{\ell}\over {1+\alpha}}\right|+S_H\right){1\over{|\alpha|^n}}
 <{{\varepsilon}\over 2},
\quad \forall x\ge 0, 
$$
\par
\noindent
whence an immediate triangle inequality gives $H(x)\to {{\ell}\over {1+\alpha}}$, from the definition of limit.
\qed

\bigskip

\par
The Theorem has the following, interesting consequence. 
\smallskip
\par 
\noindent {\bf Corollary 6.2.} ({\stampatello Converse $\SGf-$recursive formula})
\par
\noindent
{\it Let $G:\N\to\C$ be a multiplicative function and suppose that $\SG b$ converges for some $b\in\N$. Then $\SG {bc}$ converges for all $c$ coprime with $b$ and with no bad prime divisor, i.e. $\P(c)\cap \Bad = \emptyset$.
Moreover, for the same $b,c\in\N$ we have:
$$
\SG {bc} = \prod_{p|c} (1-G(p))^{-1} \cdot \SG b.
$$
}
\par
\noindent {\it Proof.} Let us first suppose $c=p$ is prime, coprime with $b$. Then the partial sums (see $\S2$ definitions) are 
$$
\SGx b x = \sum_{{r\le x}\atop {(r,pb)=1}}G(r)\mu(r)+\sum_{{r\le x}\atop {{(r,b)=1}\atop {r\equiv 0\bmod p}}}G(r)\mu(r)
= \SGx {bp} x - G(p) \SGx {bp} {x/p},
\quad 
\forall x\ge 0.
$$
\par
\noindent
Since $p$ is not a bad prime, we either have $|G(p)|<1$ or $|G(p)|>p$. Thus the hypotheses of Theorem 6.1 hold, so the convergence of $\SG b$ implies the convergence of $\SG {bp}$. Since $G(p)\neq 1$ we have more precisely $\SG {bp} = (1-G(p))^{-1}\SG b$. The result follows  by induction on the number of prime factors of $c$.  
\qed

\vfill
\eject

\par				
\noindent{\bf 7. Finiteness convergence Theorem: proof of Theorem 1.1} 
\bigskip
\par
\noindent
The implication (1)$\Rightarrow$(2) is trivial, while the implication (3)$\Rightarrow$(4) is a direct consequence of Corollary 6.2. We prove (4)$\Rightarrow$(1) in next subsection. Last but not least, (2)$\Rightarrow$(3) follows in $\S7.2$ from an important formula. 

\bigskip

\par
\noindent{\bf 7.1 Proof of the implication (4)$\Rightarrow$(1) in Theorem 1.1}
\bigskip
\par
\noindent {\it Proof of (4)$\Rightarrow$(1).} 
By assumption (4) we have that $\SG b$ converges for each $b\in\N$ coprime with the hyperbad primes of $G$. By Corollary 3.3 we deduce that $\RG b$ converges for the same $b$. In order to prove that $\RG a$ converges for all $a\in\N$ we need to deal with hyperbad prime factors. The Corollary 3.4 can do that, because $w_{p,G}=\infty$ for all $p\in \hyperB$.\hfill QED 

\bigskip

\par
\noindent{\bf 7.2 Proof of the implication (2)$\Rightarrow$(3) in Theorem 1.1}
\bigskip
\par
\noindent
The proof of the implication (2)$\Rightarrow$(3) relies on a special formula for the $\corsivoR \corsivoS \corsivoL$ series, namely $\corsivoF_G$. 
\smallskip
\par
\noindent {\bf Theorem 7.1.} ({\stampatello The $\corsivoF_G-$transformation formula}) 
{\it Let $a,b,c\in\N$ and fix a prime $p$, not dividing $abc$. Let $G:\N\to\C$ be multiplicative and $w\in\N$ be such that, say, $\Delta:=G(p^w)G(p)-G(p^{w+1})\neq 0$. Then
$$
\corsivoF_G(a,pb,c)(x)={{G(p)\corsivoF_G(a,b,p^w c)(x)-G(p^{w+1})\corsivoF_G(a,b,c)(x)}\over {\Delta}},
\quad \forall x\ge 0. 
$$
\par
\noindent
Assuming, also, that both\enspace $\corsivoF_G(a,b,c)$\enspace and\enspace $\corsivoF_G(a,b,p^w c)$\enspace converge, then\enspace $\corsivoF_G(a,pb,c)$\enspace converges, too.
}
\smallskip
\par
\noindent {\it Proof.} Since $p$ does not divide $a$, we have that $c_r(a)\neq 0$ only if $v_p(r)\leq 1$. Then for every $v\in\N_0$ the truncated $\corsivoR \corsivoS \corsivoL$ series with arguments $a,b,p^vc$ splits into two sums as follows
$$
\corsivoF_G(a,b,p^v c)(x)=\sum_{{r\leq x}\atop {(r,b)=1}}G(p^v cr)c_r(a)
= 
 \sum_{{r\leq x}\atop {(r,pb)=1}}G(p^v cr)c_r(a) 
 +
 \sum_{{r\leq x/p}\atop {(r,pb)=1}}G(p^{v+1}cr)c_{pr}(a).
$$ 
\par
\noindent
By basic properties of Ramanujan sums we have that $c_{pr}(a) = -c_r(a)$. 
Therefore specializing the above equation for $v=0$ and $v=w$ and using the multiplicativity of $G$ we get
\medskip
\quad\quad $\corsivoF_G(a,b,c)(x)\quad = \quad \corsivoF_G(a,pb,c)(x) - G(p)\corsivoF_G(a,pb,c)(x/p)$,
\medskip
\quad\quad  $\corsivoF_G(a,b,p^w c)(x) = G(p^w)\corsivoF_G(a,pb,c)(x) - G(p^{w+1})\corsivoF_G(a,pb,c)(x/p)$.
\medskip
\par
\noindent
Since 
$$
\Delta\defineq \det \left( {{1}\atop {G(p^w)}}\  {{-G(p)}\atop {-G(p^{w+1})}}\right)
=G(p^w)G(p)-G(p^{w+1})\neq 0, 
$$
\par
\noindent
we can solve for $\corsivoF_G(a,pb,c)(x)$:
$$
\corsivoF_G(a,pb,c)(x) = \left(G(p)\corsivoF_G(a,b,p^w c)(x)-G(p^{w+1})\corsivoF_G(a,b,c)(x)\right)/\Delta.
$$
\qed

\medskip

We are now ready to finish the proof of Theorem 1.1. 
\smallskip
\par
\noindent {\it Proof of (2)$\Rightarrow$(3).} 
By assumption (2) we know that $\RG a$ converges for all $a\mid N(G)$. 
By the $\corsivoL \corsivoR$ formula, Corollary 3.7, $\LG a = \corsivoF_G(1,1,a)$ converges for all $a\mid N(G)$. 
In particular, 
$\corsivoF_G(1,1,a)$ converges for all $ a$ of the form $a=\prod_{p\in J} p^\wg$, for each subset $J\subseteq \P(N(G))$ of simply bad primes of $G$. 
We now prove the following statement by induction on the cardinality of $I$: for every disjoint subsets $I,J\subseteq \P(N(G))$ we get 
$$
\psi(I,J)\defineq \corsivoF_G\left( 1, \prod_{p\in I} p, \prod_{p\in J} p^ \wg\right) 
$$
\par
\noindent
convergence. If $I=\emptyset$, this is what we just proved a few lines above. 
For the induction step, first note that 
$$
G(p^{\wg+1}) \neq G(p) G(p^\wg), 
$$
\par
\noindent
by the definition of $\wg$.  Then Theorem 7.1 tells us that the convergence of $\psi(I\cup\{p\},J)$ follows from the convergence of both $\psi(I,J)$ and $\psi(I,J\cup\{p\})$, whenever $p\notin I\cup J$. 
By induction on the cardinality of $I$ we deduce that $\psi(I,J)$ converges for all disjoint $I,J\subseteq \P(N(G))$. 
In particular 
$$
\psi (\P(N(G)),\emptyset) = \corsivoF_G (1,\rad N(G),1) = \SG {N(G)}
$$
\par
\noindent
converges.\hfill QED

Hence, the proof of Theorem 1.1 is complete.\qed

\vfill
\eject

\par				
\noindent{\bf 8. Ramanujan series factorization and Euler-Selberg products}
\bigskip
\par
\noindent
The formul\ae \enspace used in the proof of Theorem 1.1 will be used again, in their limiting case, to prove the following general formula for the values of Ramanujan expansions with multiplicative coefficients. 
\medskip
\par
\noindent{\bf 8.1 The Ramanujan series factorization and an auxiliary integers' factorization} 
\medskip
\par
\noindent
In order to state this result, we introduce, once fixed the multiplicative $G$, a kind of unique factorization for all natural numbers $a$, that we call the {\stampatello Ramanujan factorization} (or {\it $G-$Ramanujan factorization}): 
$$
a=h\cdot t\cdot\;\widetilde{a}
=h(a)t(a)\;\widetilde{}\,(a), 
$$
\par
\noindent 
where the say {\it hyperbad-factor} $h:\N\to\C$, the say {\it simply-transparent-factor} $t:\N\to\C$ and, say, the {\it regular-factor } $\;\widetilde{}\,:\N\to\C$ are the C.M.(completely multiplicative) functions defined as ($\oplus\defineq$disjoint union)
$$
h(a)\defineq \prod_{p|a,p\in \hyperB}p^{v_p(a)},
\quad 
t(a)\defineq \prod_{p|a,p\in \simplyT}p^{v_p(a)}, 
\quad
\widetilde{}\,(a)\defineq \prod_{p|a,p\not \in \hyperB \oplus \simplyT}p^{v_p(a)}. 
$$
\par
\noindent 
See that, in the following, we'll use the property: 
$$
(a,b)=1
\quad \Rightarrow \quad
(h(a),h(b))=1, \enspace 
(t(a),t(b))=1, \enspace 
(\;\widetilde{}\,(a),\;\widetilde{}\,(b))=1. 
$$
\medskip
\par
We are ready to see how this factorization of natural numbers gives a $\RGf$ factorization. 
\smallskip
\par
\noindent {\bf Theorem 8.1.} ({\stampatello Ramanujan series factorization}) 
\par
\noindent
{\it 
Let $G:\N\to\C$ be a multiplicative function such that $\RG a$ converges for all $a\in\N$. Then, writing each $a\in\N$ uniquely, following its $G-$Ramanujan  factorization, $a = h\cdot t\cdot \widetilde{a}$, we get, $\forall a\in\N$, 
$$
\RG a = 
\DG {h}
\cdot
\CGf \left( { {\rad N_T(G)} \over {\rad t} }\right) 
{ {\EG t} \over {\EG {N_T(G)}} }
\cdot
{ {\EG {\widetilde{a}}} \over {\CG {\widetilde{a}}} }  
\cdot
\RG {N_T(G)}.
$$
}
\medskip
\par
\noindent {\bf Remark 4.} {\it Notice that $\EG {N_T(G)}\neq 0$ and $\CG{\widetilde{a}}\neq 0$, so the formula makes sense. We have, in fact, 
$$
\CG {\widetilde{a}} = \prod_{p\mid \widetilde{a}}(1-G(p))
$$
\par
\noindent
which is nonzero because $G(p)\neq 1$ for all $p\mid \widetilde{a}$ (that is, no prime factor of $\widetilde{a}$ is transparent). Next, we have
$$
\EG {N_T(G)} = \prod_{p\mid N_T(G)}\sum_{k=0}^{\valg} p^k \left(G(p^k)-G(p^{k+1})\right).
$$
\par
\noindent
For every $p\mid N_T(G)$ we have $G(p^k)=G(p)^k=1$ for all $k\leq \valg$ (likewise, for $k=0$) and $G(p^{\valg+1})\neq 1$, so
$$
\EG {N_T(G)} = \prod_{p\mid N_T(G)} p^{\valg}\left(1-G(p^{\valg+1})\right)
$$  
\par
\noindent
is nonzero. 
\par
\noindent
Also, trivially, any multiplicative Ramanujan coefficient $G$ with \enspace $\RG {N_T(G)}=0$\enspace has \enspace $\RGf=\0$}.\hfill $\diamond$ 
\medskip
\par
\noindent
The Ramanujan series factorization, in particular, proves that: 
\smallskip
\par
\centerline{$G$ {\stampatello mult.Ramanujan coefficient, $\RG {N_T(G)}\neq 0$} \enspace $\Rightarrow $ \enspace $\RGf/\RG {N_T(G)}$ {\stampatello finite Euler product}.} 

\bigskip

\par
\noindent
The Proof of Ramanujan series factorization Formula simply descends from Theorem 1.1 and $\S3$ Corollaries, about factorizations of $\RGf$ and $\SGf$. 

\vfill
\eject

\par				
\noindent {\it Proof of Theorem 8.1.} Let us recall some results that we proved in the preceding sections. 
Corollary 3.4 tells us that $ \RG {bc} =  \RG b \DG c $ whenever $(b,c)=1$ and $\wg=\infty$ for all $p\mid c$. In particular we have
$$
\RG a = \DG h \RG{t\widetilde{a}}.
$$
\par
\noindent
Next, by Theorem 1.1 we have that $\SG b$ converges for any $b$ without hyperbad prime factors. In particular $\SG {t\widetilde{a}}$ converges, so Corollary 3.3 and Corollary 3.2 imply that 
$$
\RG {t\widetilde{a}} = \EG {t\widetilde{a}} \SG {t\widetilde{a}} 
\and 
\SG t = \CG {\widetilde{a}} \SG {t\widetilde{a}}.
$$
\par
\noindent
Therefore, noticing that $\SG t = \SG {\rad t}$ for trivial reasons, we get
$$
\RG a = \DG h \,\EG {t\widetilde{a}}\, (\CG {\widetilde{a}}) ^{-1} \SG {\rad t}.
$$
\par
\noindent
Since all prime factors of $t$ are simply transparent, we have $\rad t \mid \rad N_T(G)$. Then Corollary 3.2 gives
$$
\SG {\rad t} = \CGf \left( {\rad N_T(G) \over \rad t} \right) \SG {\rad N_T(G)}.
$$
\par
\noindent
By Corollary 3.3 we have $\RG {N_T(G)} = \EG {N_T(G)} \SG {\rad N_T(G)}$, whence finally the sought formula.\qed

\bigskip

\par
\noindent{\bf 8.2 The Euler-Selberg factors} 
\medskip
\par
\noindent
We prove now all the Corollaries of Theorem 8.1, needed in order to prove Corollary 1.2. 
\medskip
\par
We start from the first statement in Corollary 1.2. 
\smallskip
\par 
\noindent {\bf Corollary 8.2.} {\it Let $G$ be a multiplicative Ramanujan coefficient with transparency conductor $N_T(G)$. Then $\RG a =0$ for every $a$ that is not a multiple of $N_T(G)$.}
\smallskip
\par
\noindent
{\it Proof.} Let us write $a=h\cdot t\cdot \widetilde{a}$ as we did in Theorem 8.1. 
\par
We distinguish two cases, according to whether $\rad t$ is equal to $\rad N_T(G)$ or not. In each case, we prove that some factor in the main formula of  Theorem 8.1 vanishes.  
\par
Suppose first that $\rad t \neq \rad N_T(G) $. \enspace Then the ratio, say, \enspace $r:=(\rad N_T(G))/(\rad t)$ \enspace 
is a square-free integer $r>1$ divisible only by transparent primes. Thus 
$$
\CG  r = \prod_{p \mid r}(1-G(p)) = 0. 
$$
\par
\noindent
Using the formula in Theorem 8.1 we deduce that $\RG a = 0$ in case $\rad t \neq \rad N_T(G)$. 
\par
Let us now assume that $\rad t = \rad N_T(G)$ but $a$ is not a multiple of $N_T(G)$.  Then, there exists a prime $p$ which is simply transparent and such that 
$$
1\leq v_p(t)=v_p(a) < v_p(N_T(G)) = \valg.
$$
\par
\noindent
Since the function $\EGf$ is multiplicative and $(t,\widetilde{a})=1$, we may write
$$
\EG {t\widetilde{a}} = \EG {\widetilde{a}} \EG t = \EG {\widetilde{a}} \EG {p^{v_p(t)}} \EG {p^{-v_p(t)}t}. 
$$
\par
\noindent
The middle term in the right-most expression can be written as 
$$
\EG {p^{v_p(t)}} = \sum_{v=0}^{v_p(t)} p^v (G(p^v)-G(p^{v+1})).
$$
\par 
\noindent
Now, $p$ is a simply transparent prime, so $G(p^v)=1$ for all $v\leq\valg$. Since $v_p(t)<\valg$ by our choice of $p$, we deduce that $G(p^v)-G(p^{v+1})=0$  for each $v\leq v_p(t)$. This proves that $\EG {p^{v_p(t)}}=0$ and so $\EG {t\widetilde{a}} = 0$. 
\par
\noindent
By the formula in Theorem 8.1 we deduce that $\RG a = 0$ also in this case.
\qed

\bigskip

\par				
Then, next property is an important step in the proof of Corollary 1.2.
\smallskip
\par 
\noindent {\bf Corollary 8.3.} {\it Let $G$ be a multiplicative Ramanujan coefficient with transparency conductor $N_T(G)$. For every $a\in\N$ write the Ramanujan factorization \enspace $a=h\cdot t\cdot \widetilde{a}$ and define the function
$$
\MG a \defineq \DG h \cdot {\EG {t N_T(G)} \over \EG {N_T(G)}} \cdot {\EG{\widetilde{a}} \over \CG {\widetilde{a}}}\, .
$$
Then $\MG a $ is a multiplicative function of \enspace $a\in\N$.
}
\smallskip
\par
\noindent{\it Proof.} First of all, note that the denominators in the definition of $\MGf$ are nonzero, by Remark 1, so the function $\MGf$ is well-defined. 
Moreover, $\MGf$ is not equal to the constant zero function, because $\MG 1=1$. 
\medskip
\par
\noindent
We use now the definition and the property of Ramanujan factorization (see soon before above Theorem 8.1).
\medskip
\par
We already know that $\CGf$, $\DGf$ and $\EGf$ are multiplicative functions, therefore also the functions
$$
a\mapsto \DG {h(a)}
\and
a\mapsto {\EG {\widetilde{a}} \over \CG {\widetilde{a}}}
$$
are multiplicative.
It remains to show that the function $a\mapsto \EG {t(a) N_T(G)}/\EG{N_T(G)}$ is multiplicative. 
\par
\noindent
Let $a,b\in\N$ with $(a,b)=1$. Let $N_1$ be the largest divisor of $N_T(G)$ that is coprime with $t(b)$, and let $N_2=N_T(G)/N_1$. 
In other words, we decompose $N_T(G)=N_1\cdot N_2$ so that $(N_1,t(b))=1$ and $\rad N_2 \mid t(b)$. 
Since $t(a)$ and $t(b)$ are coprime, it follows that $(t(a),N_2)=1$. 
Then by multiplicativity of $\EGf$ we have
$$
{\EG {t(ab) N_T(G)} \over \EG {N_T(G)} } 
=
{\EG {t(a) N_1} \over \EG {N_1} } 
{\EG {t(b) N_2} \over \EG {N_2} } .
$$
\par
\noindent
Note that $\EG {N_T(G)}\neq 0$ implies $\EG {N_i}\neq 0$ for $i=1,2$ as well. 
\par
Hence, last expression is 
$$
=
{\EG {t(a) N_1} \over \EG {N_1} } 
{\EG {N_2} \over \EG {N_2} } 
{\EG {t(b) N_2} \over \EG {N_2} } 
{\EG {N_1} \over \EG {N_1} } 
=
{\EG {t(a) N_T(G)} \over \EG {N_T(G)} } 
{\EG {t(b) N_T(G)} \over \EG {N_T(G)} }.
$$
\par
\noindent
Thus the function $\MG a$ is a product of three multiplicative functions, so it's multiplicative, too.
\qed

\bigskip

\par
Another important step to prove Corollary 1.2, now. 
\smallskip
\par 
\noindent {\bf Corollary 8.4.} {\it Let $G$ be a multiplicative Ramanujan coefficient with transparency conductor $N_T(G)$. Then, there exists a multiplicative function $\MGf:\N\to\C$ such that
$$
\RG {aN_T(G)} = \RG {N_T(G)} \MG a, 
\qquad \forall a\in \N. 
$$
\par
\noindent
If \enspace $\RG {N_T(G)}\neq 0$, then the function $\MGf$ is unique, being given by the formula in Corollary 8.3.
}
\smallskip
\par
\noindent {\it Proof.} We let $\MGf$ be the function defined in Corollary 8.3 and we are going to compare it with the formula in Theorem 8.1. We apply again the Ramanujan factorization, $\forall a\in \N$. 
\par
For all natural $a$, ${\rm rad}(t(a)N_T(G)) = \rad N_T(G)$ and so 
$$
\CGf \left( { {\rad N_T(G)}\over {{\rm rad}(t(a)N_T(G))} } \right) = 1.
$$
\par
\noindent
Moreover, for every $a\in\N$, we have, from : $t$ is a C.M. arithmetic function, 
$$
t(aN_T(G)) = t(a)t(N_T(G))
= t(a)N_T(G) 
$$
\par
\noindent
and so, from :\enspace $\;\widetilde{\quad}$\enspace is C.M.,
$$
\EG {t(aN_T(G))\;\,\widetilde{}\,(aN_T(G))} = \EG {t(a)N_T(G)}\  \EG {\widetilde{a}}, 
$$ 
\par				
\noindent
by multiplicativity of $\EGf$. When we pass from $a$ to $aN_T(G)$, the not-simply-transparent factors don't change. 
\par
\noindent
It follows from the formula of Theorem 8.1 that $\RG {a N_T(G)} = \RG {N_T(G)} \MG a$. The statement about the uniqueness is clear.\qed

\bigskip

\par
Finally, we prove Corollary 1.2.
\smallskip
\par
\noindent {\it Proof.} By Corollary 8.2 we have $\RG a=0$ for every $a$ that is not a multiple of $N_T(G)$. 
\par
\noindent
Let now $a$ be a fixed natural number. From Corollary 8.4 we know that $\RG {aN_T(G)} = \RG {N_T(G)} \MG a$, where $\MGf$ is some multiplicative function; also, if $\RG {N_T(G)}\neq 0$, it is in Corollary 8.3. We now compute the values of $\MGf$, at the prime powers dividing $a$, and check that they correspond to the Euler factors displayed in Corollary 1.2. 
\par
The fact that we have two left-hand sides (first $\RG {aN_T(G)}$, then $\RG a$) in the formul\ae, for Euler products, will be a real difference only when considering the factors with simply transparent primes (as $aN_T(G)$ and $a$ have the same not-simply-transparent factors, in their Ramanujan  factorization). 
\medskip
\par
\noindent
{\it Case $p$ hypertransparent.} If $p$ is a hypertransparent prime, then $\valg =\infty$ and $G(p^v)=1$ for all $v\in\N_0$. Thus 
$$
\MG {p^{v_p(a)}} = \DG{p^{v_p(a)}} = \sum_{v=0}^{v_p(a)} p^v
$$
\par
\noindent
and this agrees, in Corollary 1.2, both with $\MGf$ definition and the corresponding factor (with $\valg=\infty$).
\medskip 
\par
\noindent
{\it Case $p$ hyperbad, not hypertransparent.} In this case $G(p)\neq 1$ and $G(p^v)= G(p)^v$, for all $v\in\N_0$. Then 
$$
\MG {p^{v_p(a)}} = \DG{p^{v_p(a)}} = \sum_{v=0}^{v_p(a)} p^v G(p^v). 
$$
\par
\noindent
This agrees, in Corollary 1.2, both with $\MGf$ definition and with the corresponding factor because $\valg=0$ (since $G(p)\neq 1$) and, of course, 
$$
G(p^v) = {G(p^v) - G(p^{v+1}) \over 1 - G(p)}.
$$
\medskip
\par
\noindent
{\it Case $p$ neither hyperbad nor transparent.} Let $p$ be a prime that is not transparent (so $G(p)\neq 1$ and $\valg=0$) and not hyperbad. Then
$$
\MG {p^{v_p(a)}} = {\EG{p^{v_p(a)}} \over \CG {p^{v_p(a)}}}  = {1 \over 1- G(p)} \sum_{v=0}^{{v_p(a)}} p^v(G(p^v) - G(p^{v+1})) 
$$
\par
\noindent
and this agrees, in Corollary 1.2, both with $\MGf$ definition and with the corresponding factor. 
\medskip
\par
\noindent
{\it Case $p$ simply transparent.} Finally, let $p$ be a simply transparent prime, so that $1\leq \valg <\infty$. 
We have 
$$
v_p(aN_T(G)) = v_p(a) + \valg, 
$$
\par
\noindent
from the complete additivity of $v_p$, because $v_p(N_T(G))=\valg$ by definition of transparency conductor. 
Now, for all natural $V$ (eventually following sums are empty, giving $0$) : 
$$
\EG {p^V} = \sum_{K=\valg}^{V} p^K\left(G(p^K)-G(p^{K+1})\right), 
$$
\par
\noindent
by the definition of $\valg$, whence 
$$
\EG {p^{v_p(a)+\valg}} = \sum_{K=\valg}^{v_p(a)+\valg} p^K\left(G(p^K)-G(p^{K+1})\right)
= \sum_{v=0}^{v_p(a)} p^{v+\valg}\left(G(p^{v+\valg})-G(p^{v+\valg+1})\right), 
$$
\par				
\noindent
and, for the same reason, 
$$
\EG {p^\valg} = p^\valg \left(1-G(p^{\valg+1})\right). 
$$
\par
\noindent
Since we may assume $\RG {N_T(G)}\neq 0$ (otherwise formul\ae's both sides vanish), applying Corollary 8.4, we get $\MG {p^{v_p(a)}}$, as defined in Corollary 8.3, is 
$$
\MG {p^{v_p(a)}} = {\EG {p^{v_p(a)+\valg}} \over \EG {p^{\valg}}} 
=\left(1-G(p^{\valg+1})\right)^{-1}\sum_{v=0}^{{v_p(a)}} p^v(G(p^v) - G(p^{v+1})). 
$$
\par
\noindent
This agrees with the definition of $\MGf$ in Corollary 1.2. 
\par
The Euler factor made of simply transparent primes is calculated, this time, distinguishing in last formula for $\RG a$, in Corollary 1.2, in two sub-cases : $a\not \equiv 0\bmod N_T(G)$ and $a\equiv 0\bmod N_T(G)$. 
\medskip
\par
\noindent
In sub-case $a\not \equiv 0\bmod N_T(G)$, the LHS $\RG a = 0$, see Corollary 8.2, implies that the RHS contains a vanishing factor; in fact, $a\not \equiv 0\bmod N_T(G)$ implies the existence of a simply transparent prime $p$, for which $v_p(a)<v_p(N_T(G)) = \valg$ : since $G(p^v) = G(p^{v+1}) = 1$ for all $v\le v_p(a)$, in this sub-case we have 
$$
\sum_{v=0}^{v_p(a)}p^v {{G(p^v)-G(p^{v+1})}\over {1 - G(p^{\valg+1})}}= 0. 
$$
\par
\noindent
In sub-case $a\equiv 0\bmod N_T(G)$, write it $a=bN_T(G)$, so that Corollary 8.4 gives 
$$
\RG a = \RG {bN_T(G)} 
= \RG {N_T(G)} \MG b, 
$$
\par
\noindent
where now $\MG b = \MG {a/N_T(G)}$ has Euler $p-$simply transparent factor (compare formul\ae, from $\EG {p^V}$ above): 
$$
{\EG {p^{v_p(b)+\valg}} \over \EG {p^{\valg}}} = {\EG {p^{v_p(a)}} \over \EG {p^{\valg}}} 
={1\over {p^\valg \left(1-G(p^{\valg+1})\right)}} \sum_{K=\valg}^{v_p(a)}p^K\left(G(p^K)-G(p^{K+1})\right)
$$
$$
=\sum_{v=0}^{v_p(a)}p^{v-\valg} {{G(p^v)-G(p^{v+1})}\over {1 - G(p^{\valg+1})}}, 
$$
\par
\noindent
appearing in second Euler product formula in Corollary 1.2.\qed

\vfill
\eject

\par				
\noindent{\bf 9. Multiplicative Ramanujan clouds}
\bigskip
\par
\noindent
We now give a description of the full multiplicative Ramanujan cloud $\cloudF_M$ of a given arithmetic function $F$. 
This is, by definition, the set of multiplicative $G\arithf$ such that $F(a) = \RG a $ for all $a\in\N$. 

\bigskip

\par
\noindent {\bf 9.1 Multiplicative Ramanujan cloud of the null function}\bigskip\noindent
\par
\noindent
We start with the important special case $F=\0$: we show that multiplicative Ramanujan coefficients $G$ of the constant null function are characterized by the convergence and vanishing of a single series that involves squarefree values of $G$.  This result extends and completes the previous work [CG] of the authors, where a description of the special multiplicative Ramanujan cloud $\cloud0_{SM}$ was given. Recall $N(G)$ in Definition B, $\S1$.
\bigskip 
\par 
\noindent {\bf Theorem 9.1 } 
{\it 
Let $G$ be a multiplicative function with Ramanujan conductor $N(G)$. We have
$$
G\in \cloud0_M 
\quad \Longleftrightarrow \quad
\sum_{(r,N(G))=1} G(r) \mu(r) = 0.
$$
}
\bigskip
\par
\noindent
{\it Proof.} Let $G$ be a multiplicative Ramanujan coefficient of $\0$ and let $N_T(G)$ be its transparency conductor. By the Finiteness Convergence Theorem we have that  $\SG {N(G)}$ and  $\SG {N_T(G)}$ converge. 
Furthermore, we have that $\EG{N_T(G)} \neq 0$ by Remark 4. Thus by the $\corsivoR \corsivoS$ and the $\corsivoS_G-$recursive formul\ae, we have 
$$ 
\SG {N(G)} = \CG{N(G)/N_T(G)} \SG {N_T(G)} = {\CG{N(G)/N_T(G)} \over \EG{N_T(G)}} \RG {N_T(G)} = 0.
$$
Since $\SG {N(G)} \defineq \sum_{(r,N(G))=1} G(r) \mu(r)$, the implication $\Rightarrow$ is proved.
\par
\noindent
Now let $G$ be a multiplicative function such that $\SG {N(G)}$ converges to zero. By the Finiteness Convergence Theorem we have that $\RG a$ converges for all $a\in\N$. Note that 
$$
\CG{N(G)/N_T(G)}  = \prod _{p\in\simplyB\setminus\simplyT} (1-G(p)) \neq 0. 
$$
Hence the previous computation gives
$$
\RG {N_T(G)} = {\EG{N_T(G)} \over \CG{N(G)/N_T(G)}} \SG {N(G)}  = 0.
$$
But then the Euler-Selberg product formula (Corollary 1.2) proves that $\RG a = 0$ for all $a\in \N$. 
\qed
\bigskip
\bigskip 
\par 
\noindent {\bf Remark 5.} {\it By the $\corsivoS_G-$recursive formula, we note that to have $\SG {N(G)}=0$ it is sufficient that the condition $\SG N=0$ holds for some arbitrary  multiple $N$ of $\rad N(G)$.\hfill $\diamond$
} 
\bigskip
\par
\noindent {\bf Remark 6.} {\it Summarizing the conclusions of  Theorem 9.1 and Remark 9.1.1, we have that each multiplicative Ramanujan coefficient $G$ of the constant null function can be constructed via the following two-step procedure: 

\medskip
 \item {(i)} choose a natural number $N$ and complex numbers  $G(p)$ for each $p$ prime, so that, defining $G$ at squarefree arguments by multiplicativity, we have $\SG N=0$;
\medskip
 \item {(ii)}  choose  the values of $G$ at all other natural arguments ensuring that $G$ is multiplicative and that  $\rad N(G)$ divides $N$.\hfill $\diamond$
} 
\medskip 

\bigskip 
\par
\noindent {\bf Remark 7.} {\it The condition  $\rad N(G)\mid N$ in step (ii) above can be restated as follows:  $\wg=\infty$ for each prime $p$ coprime with $N$ and such that $1\leq |G(p)|\leq p$. This is essentially a condition imposed on $G\in\cloud0_M$ at all powers of these (finitely many) primes. The values of $G$ at prime powers (with exponent greater than 1) at all other primes can be chosen arbitrarily.\hfill $\diamond$
} 

\vfill
\eject

\par				
\noindent{\bf 9.2 Multiplicative Ramanujan clouds of non-null functions}
\bigskip
\par
\noindent
Finally, we are going to describe  the Ramanujan coefficients of a semi-multiplicative function $F\neq \0$. 
As we shall see in the following theorems, there is an important difference with respect to the case of the null function $\0(a)\defineq 0$, $\forall a \in\N$. 

As we discovered in the previous paragraph, the values of $G\in\cloud0_M$ at prime arguments are subject to a single linear equation (namely, $\SG {N(G)}=0$), whilst the values of $G$ at higher prime powers can be essentially arbitrary (there is only an additional technical condition on the complete multiplicative index of bad primes). 
Instead, if $G\in\cloudF_M$ is a multiplicative Ramanujan coefficient of a non-null function $F\neq \0$, the Euler-Selberg product formula (Corollary 1.2) gives a recursion for the values $G(p^k)$ for fixed $p$ and varying $k$.  
This implies that the values of $G$ at high powers of primes are not arbitrary: on the contrary, they are completely determined by the values at small powers of primes. 

In order to make this observation precise, we shall use the following notion. 

\bigskip
\par
\noindent {\bf Definition 9.2} ({\stampatello opacity core})
{\it Let $G\arithf$ be a multiplicative function with finitely many transparent primes (Definition A, $\S2$) and transparency conductor equal to $N_T(G)$. 
We say that the {\stampatello opacity core} of $G$ is the function $H_G\arithf$ given by the formula 
$$
H_{G}(q) \defineq G(q\,\nt)\mu^2(q).
$$
}

\medskip
\par
\noindent {\bf Remark 8.} {\it We recall that the transparency index of $G$ at a prime $p$ is the largest $\valg \in\N\cup\{0,\infty\}$ such that $G(p^v)=1$ for each natural number $v\leq \valg$. If $p$ is not hypertransparent, we see that 
$$
H_G(p) = G(p^{v_p(\nt)+1})= G(p^{\valg+1})
\leqno{(\ast_H)}
$$ 
is the first value of $G$ at a power of $p$ that is not equal to 1. 
In fact, $H_G$ is the unique multiplicative function supported on the squarefree numbers that satifies $H_G(p)=1$ for each hypertransparent prime $p$ and that satisfies ($\ast_H$) for all other primes.\hfill $\diamond$ 
}
\medskip

It turns out that a multiplicative Ramanujan coefficient of $F\neq \0$ is completely determined by its opacity core. This is expressed qualitatively by the following rigidity statement, which we prove in section 9.3.

\bigskip\par\noindent {\bf Theorem 9.3} ({\stampatello rigidity})
{\it 
Let $F:\N\to\C$ be a non-null function and let $G_1, G_2 \in\cloudF_M$ be two multiplicative Ramanujan coefficients of $F$. 
Suppose that $G_1$ and $G_2$ have the same opacity core $H_{G_1}=H_{G_2}$.  
Then $G_1=G_2$.
}
\medskip 

\par
\noindent {\bf Remark 9.} {\it In addition to this, we shall prove in Proposition 9.6 an explicit formula that recovers all values of $G\in\cloudF_M$ from its opacity core $H_G$.\hfill $\diamond$ 
} 
\medskip
 
Thanks to Theorem 9.3, we see that  the problem of describing the full multiplicative Ramanujan cloud of $F$ is reduced to the problem of characterizing the possible opacity cores that occur. The precise description is quite technical but ultimately (see Theorem 9.5 below for details) it boils down to a single linear equation, namely
$$
\corsivoS_{H_G}(1) = \RG{\nt}.
$$
In order to state the necessary additional technical conditions we need to introduce one final notion, that of ``relative simple badness''. 

\bigskip
\par 
\noindent
{\bf Definition 9.4} ({\stampatello relative simply bad primes}) 
{\it 
Let $M,H\arithf$ be multiplicative functions and let $p$ be prime. 
We say that $p$ is a {\stampatello relative simply bad prime} for $(M,H)$ if and only if 
$1\leq |H(p)|\leq p$ and at least one of the following two conditions holds:
\medskip
\centerline {%
 (i) \qquad $w_{p,M*\mu}\neq \infty$, \qquad\qquad or %
\qquad\qquad%
 \quad (ii) \qquad  $H(p)\neq (M(p)-1)/p$.%
 }
}
\bigskip
\par				
\noindent {\bf Remark 10.} {\it In the previous definition we should think $H=H_G$ as the opacity core of some multiplicative function $G\in\cloudF_M$, and $M=M_F$ as the multiplicative function associated to the non-null semi-multiplicative function $F$ (see Proposition 1.3, $\S1$). In this case the conditions (i) and (ii) are naturally equivalent to the statement that $p$ is not a completely multiplicative prime of $G$ (see Lemma 9.7 below).\hfill $\diamond$ 
}
\medskip

\noindent 
Finally, we are able to state the classification theorem for the multiplicative Ramanujan coefficients of non-null functions.

\bigskip
\par 
\noindent
{\bf Theorem 9.5} ({\stampatello classification})
{\it 
Let $F$ be a non-null semi-multiplicative function and write $F$ in the standard form $F(a)=F(a_F) M_F(a/a_F)$ for all $a\in\N$, where $F(a_F)\neq 0$ and $M_F:\N\to\C$ is multiplicative. 
Let $H:\N\to\C$ be a multiplicative function supported on the squarefree numbers and let $N\in\N$ be a natural number divisible by $a_F$ and by all the relative simply bad primes of $(M_F,H)$. 
If 
\medskip
\centerline{$\displaystyle \sum_{(r,N)=1} H(r)\mu(r)$ converges \and  $\displaystyle  \sum_{q=1}^\infty H(q)\mu(q) = F(a_F)$,}
\medskip
\noindent
then there exists a (unique) multiplicative Ramanujan coefficient $G\in\cloudF_M$ of $F$, such that $H_G=H$. 
}
\bigskip

\par
\noindent {\bf Remark 11.} {\it An explicit formula for $G$ in terms of $H$ and $F$ is given in Propositions 9.6. We then state and prove a more precise version of Theorem 9.5 in Proposition 9.8 below.\hfill $\diamond$ 
}
\medskip

\par
\noindent {\bf Remark 12.} {\it Viceversa, we show in Lemma 9.9 that the  statement in display in Theorem 9.5 holds if $H=H_G$ is the opacity core of a multiplicative Ramanujan coefficient $G$ of $F$. Therefore Theorem 9.5 gives a complete description of the multiplicative Ramanujan cloud of $F$.\hfill $\diamond$ 
}  

\bigskip
\bigskip
\par 
\noindent
{\bf 9.3 Proofs of the rigidity theorem and of the classification }
\bigskip

We shall deduce the rigidity theorem (Theorem 9.3) from the following more precise statement, that recovers all values of $G\in\cloudF_M$ from its opacity core.

\bigskip
\par 
\noindent
{\bf Proposition 9.6} 
{\it 
Let  $F$, $a_F$ and $M_F$ be as in Theorem 9.5 and let $G\in\cloudF_M$. Then for every prime $p$ and all $v\in\N$ we have the explicit formula
$$
G(p^{\vpa+v}) = H_G(p) + (1-H_G(p))\gmf(p^v)
\leqno{(\ast_G)}
$$
where $H_G$ is the opacity core of $G$ and $\gmf$ is the canonical Ramanujan coefficient of $M_F$ (see Definition E, $\S2$).
}
\bigskip

\noindent
{\it Proof.} 
If $p$ is hypertransparent, then $H_G(p)$ and $G(p^k)=1$ for all $k\in\N$, so the formula ($\ast_G$) holds. 
If $p$ is not hypertransparent, then $\valg=\vpa$ because by the Euler-Selberg product formula (Corollary 1.2) we know that $a_F$ is the transparency conductor of $G$. 
Therefore ($\ast_H$) tells us that 
$$
H_G(p) = G(p^{\vpa+1})=G(p^{\valg+1})\neq 1.
\qquad\qquad(eq:hp2)
$$ 
A simple property of the canonical Ramanujan coefficient is that $\gmf(p)=0$ for each prime $p$, therefore ($\ast_H$) is equivalent to the case $v=1$ of ($\ast_G$). 
To prove the formula for $v\geq 1$ we use the Euler-Selberg product formula (Corollary 1.2) of $F=\corsivoR_G$. 
An inspection of this factorization formula reveals that
$$
 M_F(p^v) - M_F(p^{v-1}) = p^v {G(p^{\vpa + v}) - G(p^{\vpa+v+1}) \over 1 - H_G(p)}
$$
for all $v\in\N$. We deduce that 
$$
G (p^{\vpa + v+1}) =   G(p^{\vpa+v}) - (1-H(p)) p^{-v} M_F'(p^v)
$$
\par				
\noindent
for all $v\in\N$, where $M_F'$ denotes the Eratosthenes transform (Definition C, $\S2$) of $M_F$. Since the canonical Ramanujan coefficient $\gmf$ of $M_F$ is given explicitly by
$$ \gmf(p^v) = -\sum_{k=1}^{v-1} {M_F'(p^v) \over p^v},$$
the proposition is proved.
\qed
\medskip
\par
\noindent
The rigidity theorem is now a corollary:
\smallskip
\par
\noindent
{\it Proof of Theorem 9.3.}  
From Proposition 9.6 we get  $G_1(p^v)=G_2(p^v)$ for every $v>\vpa$. On the other hand, we have that $G_1(p^v)=G_2(p^v)=1$ for each $v=1,\ldots,\vpa$. Therefore $G_1$ and $G_2$ have the same value at all powers of primes; since they are multiplicative, they coincide.
\qed
\medskip
Before we embark the proof of the classification theorem (Theorem 9.5), we state and prove an easy lemma that illustrates the rationale behind the notion of relative simply bad primes (Definition 9.4). 
Briefly, we know from Proposition 9.6 that a multiplicative Ramanujan coefficient $G$ can be fully recovered from $F=\RGf$ and $H=H_G$. 
This means that all the relevant set of primes (transparent, hypertransparent, bad and hyperbad) attached to $G$ should admit a direct description in terms of $F$ and $H$ only. By the structure theorem, the simply transparent primes of $G$ are recovered as the prime divisors of $a_F\defineq \min\{ a\in\N:F(a)\neq 0\}$. Moreover Remark 9.2.1 tells us that the hypertransparent primes $p$ are characterized by the property $H(p)=1$. If $p$ is a prime that does not divide $a_F$, we have $H(p)=G(p)$, so badness is detected by the inequality $1\leq | {H(p)} |\leq p$. It remains only to differentiate between simply bad and hyperbad primes. This is accomplished by Definition 9.4 (see also Remark 9.4.1) and by the following lemma. 
\medskip
\par 
\noindent
{\bf Lemma 9.7} 
{\it 
Let $G,M\arithf$ be multiplicative functions and let $G_M$ be the canonical Ramanujan coefficient of $M$. 
Let $p$ be some prime number such that $G(p)\neq 1$ and suppose that 
$$
G(p^v) = G(p) + (1-G(p)) G_M(p^v)
$$
for all $v\in\N$. 
Then the following are equivalent: 
\medskip
\quad (i)\ \  \quad  $\wg=\infty$; 
\medskip
\quad (ii) \quad $w_{p,M*\mu} = \infty$ \quad  and\quad   $G(p) = (M(p) - 1 )/ p$.
}
\bigskip 
\par
\noindent
{\it Proof.} 
First we recall that the values of the canonical Ramanujan coefficient at prime powers of $p$ are 
$$
G_M(p^v) = -\sum_{k=1}^{v-1} {M'(p^v) \over p^v},
$$
where $M'\defineq M\ast \mu$ is the Eratosthenes transform of $M$. Then for all $v\in\N$ we have 
$$
 {G(p^{v}) - G(p^{v+1}) \over 1 - G(p)}
 = { M'(p^v) \over p^v}.
$$ 
If (i) holds, then $G(p^v)=G(p)^v$ for all $v\in\N$. 
Therefore
$$
{M'(p^{v})  }
= 
 p^v {G(p)^{v} - G(p)^{v+1} \over 1 - G(p)} = p^v G(p)^v
 $$
for all $v\in \N$. 
This means that $M'(p)= p G(p)$ and that $M'$ is completely multiplicative along $p$. Since $M'(p)= M(p)-1$, we get (ii).

\noindent 
Suppose now that (ii) holds. 
Then 
$
M'(p^{v}) = (M'(p))^v = p^v G(p)^v
$
for all $v\in\N$ and so 
$$
 G(p^{v}) = G(p)  - (1-G(p))\sum_{k=1}^{v-1} G(p)^v.
$$
\par				
\noindent
Using the formula for the sum of a geometric progression we verify that this simplifies to $G(p^v)=G(p)^v$, and therefore we get (i).
\qed
\bigskip

We now present the proof of the classification theorem, in the following more precise form:

\bigskip
\par 
\noindent
{\bf Proposition 9.8} {\it Let  $F$, $a_F$, $M_F$ and $H$ be as in Theorem 9.5. 
Let $\gmf$ be the canonical Ramanujan coefficient of $M_F$ and let $G\arithf$ be the unique multiplicative function such that
$$ G(p^{v+v_p(a_F)}) = H(p) + (1- H(p)) \gmf(p^{v}) $$
for each $p$ prime and all $v\in\N$, and such that $G(p^v)=1$ for all prime $p\mid a_F$ and all integers $v\in [1,v_p(a_F)]$. 
Then $H$ is the opacity core of $G$ and $G\in\cloudF_M$. 
}

\bigskip 
\noindent
{\it Proof.} 
We divide the proof in 7 steps.

\medskip
\par 
\noindent
{\bf Step 1: $H=H_G$.} 
We recall that the canonical Ramanujan coefficient satisfies $\gmf(p)=0$ for every prime $p$, so putting $v=1$ in the definition of $G$ we obtain
$$
H(p) = G(p^{\vpa +1}).
$$
If $H(p)\neq 1$ then $H(p)$ is the first value of $G$ at a power of $p$ that is not equal to $1$. 
If $H(p)=1$ instead, then $G(p^v)=1$ for all $v\in\N$. 
Since moreover $H$ is multiplicative and supported on squarefree numbers, we deduce that $H$ is the opacity core of $G$ (see Remark 9.2.1).

\medskip
\par 
\noindent
{\bf Step 2: $a_F=N_T(G)$.}  
By hypothesis we have that $\SH(N)$ converges and that $\SH(1)= F(a_F)\neq 0$. 
By the $\corsivoS_G-$recursive formula we also know that
$$
\SH(N)\prod_{p\mid N} (1-H(p)) = \SH(1),
$$
hence $H(p)\neq 1$ for all prime divisor $p$ of $N$. 
In particular, for every $p\mid a_F$ we have $G(p^{\vpa+1})=H(p)\neq 1$ and $G(p^v)=1$ for all $v\leq \vpa$. 
This means that for each $p\mid a_F$ we have
$$
v_p(\nt)=\valg=\vpa.
$$
If instead $p$ does not divide $a_F$, then $G(p)=H(p)$ and either $\valg=0$ (if $H(p)\neq 1$) or $\valg=\infty$ (when $H(p)=1$). 
In both these cases we get $v_p(\nt)=\vpa=0$. 
Therefore $v_p(\nt)=\vpa$ holds for every prime $p$, which implies $a_F=\nt$.

\medskip
\par 
\noindent
{\bf Step 3: $\rad N(G)\mid N$.}  
As we have already seen, all simply transparent primes of $G$ are divisors of $a_F$, and so they divide $N$. 
If $p$ is a simply bad prime of $G$ that is not simply transparent, then $\vpa=0$ and $G(p)=H(p)\neq 1$. 
Therefore we may apply Lemma 9.7 with $M=M_F$. 
Since $p$ is simply bad, we have that $\wg \neq\infty$ and so Lemma 9.7 implies that $p$ is a {\it relative simply bad prime} (see Definition 9.4) for the pair $(M_F,H)$. 
We therefore conclude that $p\mid N$ for every simply bad prime  $p$ of $G$.

\medskip
\par 
\noindent
{\bf Step 4: $\RG a$ converges $\forall a\in\N$.}  
For every prime $p$ coprime with $a_F$ we have $H(p) = G(p)$, therefore for each $b$ multiple of $a_F$ the series $\SG b$ coincides with $\corsivoS_H(b)$. 
By assumption we know that $\corsivoS_H(N)$ converges, so by the $\corsivoS_G-$recursive formula we have that $\corsivoS_H(d)$ converges for all $ d\mid N$. 
Since $N$ is a multiple of $\rad N(G)$, and $N(G)$ is a multiple of $a_F=N_T(G)$, we deduce that 
$$
\SG {N(G)} = \corsivoS_H(N(G)) = \corsivoS_H(\rad N(G))
$$ converges. 
The convergence of $\RG a$ for all $a\in\N$ follows from  the finiteness convergence theorem (Theorem 1.1).

\medskip
\par 
\noindent
{\bf Step 5: $F(a_F) = \RG {N_T(G)}$.}   
Since $a_F=\nt$, we note the following equality:
$$
\EG {N_T(G)} 
= \prod_{p\mid N_T(G)} (1-G(p^{\vpa+1})) 
= \prod_{p\mid N_T(G)} (1-H(p))
= \corsivoC_H (N_T(G)).
$$
\par				
\noindent
By the $\corsivoR \corsivoS-$formula and by the $\corsivoS_G-$recursive formula we then obtain
$$
\RG {\nt} 
= \EG {\nt} \, \SG {\nt}
= \corsivoC_H(\nt)\,  \corsivoS_H(\nt)
= \corsivoS_H(1).$$
Moreover, $\corsivoS_H(1)$ is equal to $F(a_T)$ by the assumption on $H$.

\medskip
\par 
\noindent
{\bf Step 6: $\valg = \infty$ $\Rightarrow$ $M'_F(p^v)=p^v$ for all $v\in\N$.}  
We already noted (see Steps 1 and 2) that for all hypertransparent prime $p$ we have $H(p)=1$ and that $p$ does not divide $N$. 
In particular, such $p$ is not a relative simply bad prime of $(M_F,H)$. In other words, we have
$$
w_{p,M'_F}=\infty
\qquad{\stampatello and}\qquad
H(p)= (M_F(p)-1)/p.
$$ 
We recall that $w_{p,M'_F}=\infty$ means that $M'_F$ is completely multiplicative along $p$. 
Since $H(p)=1$ and $M'_F(p)=M_F(p)-1$, we conclude that $M'_F(p^v)=(M'_F(p))^v=p^v$ for all $ v\in\N$.

\medskip
\par 
\noindent
{\bf Step 7: $F(a) = \RG a$ for all $a\in\N$.}  
By the Euler-Selberg product formula (Corollary 1.2) we have that
$$
\RG a = 
\RG {\nt}
\prod_{\valg=\infty} \left(\sum_{v=0}^{v_p(a)} p^v \right) 
\prod_{\valg \neq \infty} \left(  \sum_{v=0} ^ {v_p(a)-\valg} p^{v}{G(p^{\valg+v}) - G(p^{\valg + v+1}) \over 1 - G(p^{\valg+1}}
\right)
$$ 
for all $a\in\N$. 
We now aim to compare this formula with $F(a_T)M_F(a/a_T)$. 
First note that $\RG{\nt} = F(a_F)$ by Step 5.
Next, for every prime $p$ with $\valg= \infty$ we have 
$$
\sum_{v=0}^{\vpa} p^v = 
\sum_{v=0}^{\vpa} M'_F(p^v) = M_F(p^{\vpa}) 
$$
by Step 6 and the formula $M_F=M'_F\ast 1$. 
Finally, when  $\valg\neq \infty$  we have 
$G(p^{\valg+1})=H(p)$, so 
$$
\sum_{v=0} ^ {v_p(a)-\valg} p^{v}{G(p^{\valg+v}) - G(p^{\valg + v+1}) \over 1 - G(p^{\valg+1})}
 =
\sum_{v=0} ^ {v_p(a)-\valg} p^{v}\left({\gmf(p^{v}) - \gmf(p^{v+1}) }\right)
= 
\sum_{v=0} ^ {v_p(a)-\valg} M_F'(p^v)
$$
which is equal to $M_F(p^{v_p(a)-\vpa}) $ since $\valg=\vpa$. 
Thus the Euler-Selberg product formula (Corollary 1.2) of $\RGf$ coincides with the canonical Selberg factorization of $F$.
\qed
\bigskip 

\par
\noindent
Reasoning like in Steps 4 and 5  of the above proof we get the following result.

\bigskip
\par 
\noindent
{\bf Lemma 9.9} 
{\it 
Let $G\arithf$ be a multiplicative Ramanujan coefficient with Ramanujan conductor $N(G)$, transparency conductor $N_T(G)$ and opacity core $H_G$. Then 
$$\SG{N(G)} = \corsivoS_{H_G}(N(G)) 
\qquad{\stampatello and} \qquad
\RG {N_T(G)} = \corsivoS_{H_G}(1).
$$
}

Thus we see that the series in display in Theorem 9.5 correspond to the notable series $\SG{N(G)}$ and $\RG {N_T(G)}$ attached to $G$. 

\vfill
\eject

\par				
\noindent{\bf Appendix. Odds \& ends (more details, alternative proofs, generalizations)}
\bigskip
\par
\noindent
In order to give, in next $\S A.2$, some examples of not-converging coprime series, we supply, in next $\S A.1$, a short parade of elementary analytic results for the square-free numbers (as our $\S A.2$ counterexamples need it). 

\bigskip

\par
\noindent{\bf A.1. Estimates for square-free numbers}
\bigskip
\par
\noindent
We start with a classic: the asymptotic estimate for the number of square-free numbers up to $x$, as $x\to \infty$. 
\smallskip
\par
\noindent {\bf Lemma A.1.1.} ({\stampatello square-free numbers asymptotic count}) 
\par
\noindent
{\it
For all large $x$ we have
$$
\sum_{q\leq x} \mu^2(q) = {x \over \zeta(2)} + O(\epsilon(x) \sqrt x),
$$
\par
\noindent
where $\epsilon(x) = \exp \left(- c (\log x)^{3/5} (\log \log x )^{-1/5}\right)$ for some constant $c>0$.
}
\par
\noindent {\it Proof.} This estimate was proved by Walfisz [Wa] using the classical
zero-free region estimates due to Vinogradov and Korobov through
exponential sums. (This result is quoted also in the survey paper of Pappalardi [Pa] but with a mistake : $1/5$ instead of $-1/5$).
\qed
\medskip
\par
\noindent {\bf Lemma A.1.2.} ({\stampatello Dirichlet series for square-free integers})
\par
\noindent
{\it For all $s\in\C$ with $\Re s>1$ we have 
$$
\sum_{q=1}^\infty \mu^2(q) q^{-s} = \prod_{p\in\P} (1+p^{-s}) = {\zeta(s) \over \zeta(2s)}.
$$
}
\par
\noindent {\it Proof.} This is the half plane of absolute convergence, so everything works fine. \qed
\medskip
\par
\noindent {\bf Lemma A.1.3.} ({\stampatello truncated Dirichlet series for square-free})
\par
\noindent
{\it For all $s\neq 1$ with $\Re s\geq 1/2$ we have
$$
\sum_{q\leq x} \mu^2(q) q^{-s} = {6\over \pi^2} {x^{1-s}\over 1-s} + {\zeta(s) \over \zeta(2s)} + o(1).
$$
\par
\noindent
as $x\to\infty$. 
}
\par
\noindent {\it Proof.} 
First use summation by parts, and Cauchy's criterion of convergence to have the main term and a constant dependent  on $s$. Then as $x$ goes to infinity prove that the convergence is uniform on compact sets, which implies that the constant is holomorphic in $s$. Then go to the plane of absolute convergence to see that the constant coincides with $\zeta(s)/\zeta(2s)$. Another approach is to use Perron's formula.
\qed
\medskip
\par
\noindent {\bf Lemma A.1.4.} ({\stampatello truncated Dirichlet series for square-free with coprimality conditions})
\par
\noindent
{\it Let $b\in\N$. Then for all $s\neq 1$ with $\Re s\geq 1/2$ we have
$$
\sum_{{q\leq x}\atop {(q,b)=1}} \mu^2(q) q^{-s} = C_1(b) {x^{1-s}\over 1-s} +  C_2(b) {\zeta(s) \over \zeta(2s)}+ o(1).
$$
as $x\to\infty$, where \enspace $C_1(b)\defineq {6 \pi^{-2}}\prod_{p\mid b}( 1+p^{-1})^{-1}$ \enspace and\enspace  $C_2(b)\defineq \prod_{p\mid b}( 1+p^{-s})^{-1} $. 
}
\par
\noindent {\it Proof.} Use the Lemma A.1.3 together with the abstract lemmas for convergence of truncated Euler factors. More precisely, for every square-free $b\in\N$ we define the following function defined on the positive real numbers:
$$
H_b(x) \defineq \sum_{{q\leq x}\atop {(q,b)=1}}\mu^2(q) q^{-s}  - C(b) {x^{1-s}\over 1-s}.
$$
\par				
\noindent
This consists of two parts, let $H'_b(x)\defineq \sum_{{q\leq x}\atop {(q,b)=1}} $ and $H''_b(x)\defineq (6C(b) /\pi^2)\cdot x^{1-s}/(1-s) $. 
By the truncated partial Euler product formula for multiplicative functions we have
$$
H'_1(x) = \sum_{d\mid b} A(d) H'_b(x/d),
$$
where $A(d) = \mu^2(d) d^{-s}$. Interestingly, we also have the same formula for $H''_b(x)$, for example for $b=p$ prime we have
$$
H''_1(x) = C(1) {x^{1-s}\over 1-s}  = {C(1)\over 1+p^{-1}} {{x^{1-s}}\over {1-s}} + p^{-s} {C(1)\over 1+p^{-1}} {{(x/p)^{1-s}}\over {1-s}},
$$
which is $H''_p(x) + A(p) H''_p(x/p)$. 
The more general formula $H'_1(x) = \sum_{d\mid b} A(d) H'_b(x/d)$ for every $b$ coprime is obtained similarly. 
Putting all together, we get
$$
H_1(x) = \sum_{d\mid b} A(d) H_b(x/d)
$$
\par
\noindent
for every  $b\in\N$ square-free. However, we already know by the previous lemma that $\limx H_1(x)= \zeta(s)/\zeta(2s)$. By the iterated converse convergence theorem for linear combinations of multiplicative shifts, it follows that also $\limx H_b(x)$, exists, and it is equal to
$$
\limx H_b(x) = \prod_{p\mid b}{1 \over 1+p^{-s}} \cdot {\zeta(s)/\zeta(2s)}.
$$ 
This is exactly what we wanted to prove. 
\qed

\bigskip

\par
\noindent{\bf A.2. A class of coprime series with convergence problems}
\bigskip
\par
\noindent
We give a first result, in the direction of the necessity of assuming some convergence hypotheses, in our Classification Theorem in [CG]. It gives a class of multiplicative $G:\N \rightarrow \C$ whose coprime series have convergence problems. 
\smallskip
\par
\noindent {\bf Proposition A.2.1.} ({\stampatello necessity of convergence hypotheses in Classification theorem})
\par
\noindent
{\it Let $s$ be a complex number with  $1/2\leq \Re s<1$ and fix two distinct primes $p_1,p_2\in\P$. 
Let $G$ be a multiplicative function such that 
$$
G(p_1)=p_1^{1-s},
\quad\quad\quad
G(p_2)\neq p_2^{1-s}
$$
\par
\noindent
and $G(p)=-p^{-s}$ for every $p\in\P$ with $p\neq p_1$ and $p\neq p_2$. Then for every $b\in\N$ with $(b,p_1)=1$ we have
$$
\sum_{(q,b)=1} G(q)\mu(q) \ \ {\rm\bf converges},
\quad \quad{\rm but } \quad \quad
\sum_{(r,p_1b)=1} G(q)\mu(q) = \infty.
$$
\par
\noindent
If moreover $G(p_1^2)=0$, then for every $b\in\N$ we also have that 
$$
\sum_{(r,b)=1}G(r)c_r(p_1) = \infty.
$$
}
\par
\noindent {\it Proof.} Recall $\SGx b x$ definition in $\S2$. First of all, for every $m$ coprime with $p_1p_2$ we have that $\SGx {p_1p_2m} x$ is just a truncated Dirichlet series for square-free numbers with coprimality conditions. Thus, by previous $\S A.1$ results, 
$$
\SGx {p_1p_2m} x = a_m(s) {x^{1-s}} + F_{m}(s) + o_m(1),
$$
\par
\noindent
where $F_{m}(s)$ and $a_{m}(s)$ are constant with respect to $x$ and more precisely
$$
a_{m}(s) \defineq {C_1(p_2p_2m) \over 1-s}
\quad\quad{\rm and} \quad\quad
F_{m}(s) = C_2(p_1p_2m) {\zeta(s)\over \zeta(2s)}. 
$$
\par				
\noindent
This immediately implies that the series $\SG {p_1p_2m}$ diverges. In order to deal with the other series we are going to use the truncated partial Euler factorization formula for S-series. First we finish up establishing the divergence over the numbers coprime with $p_1m$, where again $m$ is coprime with $p_1p_2$:
$$
\SGx {p_1m} x = \SGx {p_1p_2m} x - G(p_2) \SGx {p_1p_2m} {x/p_2} \enspace \sim a_m(s) x^{1-s} \left( 1 - {G(p_2)\over p_2^{1-s}}\right). 
$$
\par
\noindent
So it diverges. Now we deal with the convergence in the cases of the form $b=p_2m$:
$$
\SGx {p_2m} x = \SGx {p_1p_2m} x - G(p_1) \SGx {p_1p_2m} {x/p_1} \enspace = F_{m}(s) - G(p_1) F_{m}(s) + o(1).
$$
\par
\noindent
Then we establish the convergence of $\SG m$ when $m$ is coprime with $p_1p_2$:
$$
\SGx m x = F_{m}(s) (1-p_1^{1-s}) (1-G(p_2)).
$$
\par
\noindent
Finally, let us deal with the last assertion. First, for $b$ multiple of $p_1$ we have 
$$
\sum_{(r,b)=1} G(r) c_r(p_1) = \sum_{(r,b)=1} G(r) \mu(r)
$$ 
\par
\noindent
and this diverges by the previous results in this Proposition. More precisely
$$
\corsivoR_{p_1m}(p_1,x) \defineq \sum_{{(r,p_1m)=1\atop r\leq x}} G(r) c_r(p_1) = 
\alpha_m x^{1-s} + O(1)
$$
\par
\noindent
for some nonzero $\alpha_m$ dependent on $m$, where $m$ is any number coprime with $p_1$.  
To finish, we need to consider the case $b$ coprime with $p_1$. Since $G(p_1^2)=0$, 
$$
\corsivoR_{b}(p_1,x) = \corsivoR_{p_1b}(p_1,x) + G(p_1)c_{p_1}(p_1) \corsivoR_{p_1m}(p_1,x/p_1). 
$$
We recall that $G(p_1)c_{p_1}(p_1) = (p-1) p^{1-s}>0$. Then, the following diverges: 
$$
\corsivoR_{b}(p_1,x) \sim \alpha_m x^{1-s} (1+ (p-1) p^{1-s}) 
$$
\qed

\bigskip

\par 
\noindent{\bf A.3. A longer Theorem 1.1: two equivalent properties more}
\bigskip
\par
\noindent
We remark that our Theorem 1.1 can be completed, with other two equivalent properties, namely
\medskip
\quad (5) \enspace $\LG b$ converges for all $b$ not divisible by any hyperbad prime $p\in\hyperB$; 
\smallskip
\quad (6) \enspace $\LG a$ converges for all $a\in\N$. 
\medskip
\par
\noindent
We prove now the implications (4)$\Rightarrow$(5), (5)$\Rightarrow$(6) and (6)$\Rightarrow$(1) and \lq \lq QED\rq \rq, here, separates each proof. 
\medskip
\par
\noindent {\it Proof of (4)$\Rightarrow$(5).} 
By assumption $(4)$ we have that $\SG b$ converges for each $b\in\N$ coprime with the hyperbad primes of $G$ : by Corollary 3.5 we deduce that $\LG b$ converges for the same $b$, namely $(5)$.\hfill QED 
\smallskip
\par
\noindent {\it Proof of (5)$\Rightarrow$(6).} 
In order to prove $(6)$, it suffices to apply Corollary 3.6 to the $b\in\N$ without hyperbad primes (for these $b$, $(5)$ ensures $\LG b$ convergence) and to the $c\in\N$ made only of hyperbad primes.\hfill QED
\smallskip
\par
\noindent {\it Proof of (6)$\Rightarrow$(1).} 
We already know that $G$ is a Ramanujan coefficient iff it is a Lucht coefficient; this is given by Corollary 3.7, actually proving, in particular, that $(6)$ implies $(1)$.\hfill QED
\medskip
See that however Professor L.G. Lucht has full credit for having already proved (at least, 25 years ago), see [Lu1], the equivalence between the convergence of His series and Ramanujan series, namely, in Theorem 1.1 this \enspace (1)$\Leftrightarrow$(6) here. 
\medskip
\par
The same result is quoted also in [La], of course. By the way, the more flexible result on divisor-closed sets (implied by Corollary 3.7), maybe, could be of interest, as [La] needs, in the convergence of Ramanujan expansions in subsets of natural numbers. 

\bigskip

\vfill
\eject

\par				
\noindent{\bf A.4. A shorter, alternative proof of Corollary 3.4}
\bigskip
\par
\noindent
We sketch an alternative, slightly more direct proof of Corollary 3.4.
\smallskip
\par
\noindent {\it Proof.} Given $n,m\in\N$, we define (assuming, here, the multiplicative $G:\N \rightarrow \C$ is implicit, in notation) 
$$
\corsivoR_{G,n}(m,x)=\corsivoR_{n}(m,x)\defineq \sum_{{r\leq x\atop (r,n)=1 }} G(r) c_r(m). 
$$
\par
\noindent
Then we have $\corsivoR_{c}(a,x)=\corsivoR_{c}(b,x)$ and the following formul\ae: (with \thinspace $c_t(b)=\mu(t)$, in next second) 
$$
\RGx a x = \sum_{d\mid c \rad c} G(d) c_d(a) \corsivoR_{c}(a,x/d)
\and
\RGx b x = \sum_{t\mid c \rad c} G(t) c_t(b) \corsivoR_{c}(b,x/t).
$$
\par
\noindent
Writing $d=ht$ and using Kluyver's formula ($\S2$), like in above original proof, this time for $c_d(a)$, it is now easy to show the relation between $\RG a$ and $\RG b$ stated in Corollary 3.4. 
\qed 

\vfill
\eject

\par				
\noindent{\bf A.5. A general approach to our main argument: $M-$expansions}
\bigskip
\par
\noindent
We give a kind of general approach to the arguments we used to study Ramanujan expansions, modeling on them the following, more general (but, trying to catch main features) with respect to: $F=\RGf$. We call it an {\stampatello $M-$expansion}, say, $F(a)=\sum_{q=1}^{\infty}M_q(a)$, where the $M_q(a)$ is called an {\stampatello $M-$term}, by definition, when it is an arithmetic function of two arguments ($a,q\in\N$), satisfying the axioms:
\bigskip
\par
\item{(0)} ({\stampatello multiplicativity}) : $M_q(a)$ is multiplicative w.r.t. $q\in \N$; 
\medskip
\item{(1)} ({\stampatello indipendence}) : $\forall a\in\N$, $\forall p\in\P$, $\forall k\in\N_0$, $M_{p^k}(a)=M_{p^k}(p^{v_p(a)})$ ($\Rightarrow $ $M_q(a)=M_q(1)$, $\forall (a,q)=1$); 
\medskip
\item{(2)} ({\stampatello vertical limit}) : $\forall a\in\N$, $\forall p\in\P$, $\exists K_p(a)\in \N$ : $M_{p^k}(a)=0$, $\forall k>K_p(a)$. 
\medskip
\par
\noindent
Here $V_p(a)\defineq \min\{K_p(a)\enspace \hbox{\it in}\enspace (2)\}$ is the {\bf vertical limit of $q=p^k$ with respect to $a$}, esp., for Ramanujan expansions, $V_p(a)=v_p(a)+1$: this is the {\stampatello Ramanujan vertical limit}. We used, implicitly, this in all of this paper's results. (The results of following $\S A.6$ depend strongly on it.) From $(1)$ it is clear that, in general, $V_p(a)$ depends upon $v_p(a)$. 
\medskip
\par
With these three axioms, in this order, we get, formally, the Euler product
$$
\sum_{q=1}^{\infty}M_q(a)=\prod_p \sum_{K=0}^{\infty}M_{p^K}(a)
=\prod_p \sum_{K=0}^{\infty}M_{p^K}(p^{v_p(a)})
=\prod_p \sum_{K=0}^{V_p(a)}M_{p^K}(p^{v_p(a)})
$$
\par
\noindent
and we usually separate
$$
=\prod_{p|a}\sum_{K=0}^{V_p(a)}M_{p^K}(p^{v_p(a)})\cdot \prod_{p\nondivide a}\left(M_1(1)+M_p(1)+\ldots+M_{p^{V_p(a)}}(1)\right), 
$$
\par
\noindent
but we wish to have, in RHS, as infinite product (over $p\nondivide a$), say : \enspace $\sum_{(r,a)=1}\mu^2(r)M_r(1)$. 
\medskip
\par
\noindent
We may add other two axioms with which the {\it $M-$term} $M_q(a)$ becomes a {\stampatello simple $M-$term}, i.e. 
\medskip
\par
\item{(3)} ({\stampatello square-free support}) : $M_r(1)=\mu^2(r)M_r(1)$, $\forall r\in\N$; 
\medskip
\item{(4)} ({\stampatello normalization}) : $M_1(1)=1$.
\medskip
\par
\noindent
With these two more, previous formal calculation (esp., under absolute convergence hypothesis, see above) becomes 
$$
\sum_{q=1}^{\infty}M_q(a)=\prod_{p|a}\sum_{K=0}^{V_p(a)}M_{p^K}(p^{v_p(a)})\cdot \prod_{p\nondivide a}\left(1+M_p(1)\right)
=\prod_{p|a}\sum_{K=0}^{V_p(a)}M_{p^K}(p^{v_p(a)})\cdot \sum_{(r,a)=1}\mu^2(r)M_r(1), 
$$
\par
\noindent
which gets closer, to our calculations for Ramanujan expansions. 
\par
In case of convergence in $a\in\N$, let's define the {\stampatello $M-$series} at $a$ (no confusion should arise with the notation in Corollary 1.2), generalizing $\RG a$, like the following {\stampatello copriMe series} $\corsivoS(a)$ generalizes $\SG a$ (the above finite factor is the Euler$-M-$factor, but the convergence problems are with these two series): 
$$
\corsivoM(a)\defineq \sum_{q=1}^{\infty}M_q(a),
\qquad
\corsivoS(a)\defineq \sum_{(r,a)=1}\mu^2(r)M_r(1). 
$$
\par
\noindent
Notice that $\corsivoM(a)$, resp., $\corsivoS(a)$, reduces to $\RG a$, resp., $\SG a$, when $M_q(a)=G(q)c_q(a)$. 
\par
Since the \lq \lq easy\rq \rq, say, properties (like in $\S3$) should generalize at once, do these 5 axioms prove any of our main results (see $\S1$) ? If not, what are the necessary axioms for that aim and what are sufficient ? We don't have time (though margin, yes, it's not too narrow) to answer. 

\vfill
\eject

\par				
\noindent{\bf A.6. Finite Ramanujan expansions, purity and Hildebrand's coefficients}
\bigskip
\par
\noindent
As given, in the book [ScSp] (see Theorem 1.1 of V.1), from Hildebrand's 1984 paper [Hi], each arithmetic function $F$ has the following, say, {\stampatello Hildebrand-Ramanujan expansion} 
$$
F(a)=\sum_{q|a}\Hildebrand(q\rad q)c_{q\rad q}(a),
\qquad
\forall a\in\N.
$$
\par
\noindent
(Set, in quoted result, their $r^*=r\prod_{p|r}p=r\rad r:=q\rad q$, changing letter from $r$ to $q$, in our notation.)
\par
\noindent
We call, in fact, this $\Hildebrand:\N\to\C$ the {\stampatello Hildebrand coefficient} of our $F$. It can be given explicitly, in a recursive manner [ScSp], as (recall $1\rad 1=1$) 
$$
\Hildebrand(1)\defineq F(1),
\qquad
\Hildebrand(q\rad q)\defineq {1\over {c_{q\rad q}(q)}}\left(F(q)-\sum_{{d|q}\atop {d<q}}\Hildebrand(d\rad d)c_{d\rad d}(q)\right),\enspace \forall q>1. 
$$
\par
\noindent
Before to proceed, see that previous denominator doesn't vanish, as from H\"older's formula (see $\S2$): 
$$
c_{q\rad q}(q)=\varphi(q\rad q){{\mu(\rad q)}\over {\varphi(\rad q)}}\neq 0. 
$$
\par
\noindent
By the product formula [T] for $\varphi$, i.e., 
$$
\varphi(n)=n\prod_{p|n}\left(1-{1\over p}\right), 
$$
\par
\noindent
in fact, 
$$
\varphi(q\rad q)/\varphi(\rad q)=q, 
$$
\par
\noindent
whence : 
$$
c_{q\rad q}(q)=q\mu(\rad q), 
$$
\par
\noindent
not vanishing, from $\rad q$ being square-free (by definition). \enspace Notice, also, that ${\rm Hi}_{\0}=\0$. 
\par
How can be proved Hildebrand's (finite Ramanujan expansion) formula? Well, we'll recover it, from the vertical limit of Ramanujan sums: 
$$
c_{p^K}(p^{v_p(a)})\neq 0
\quad \Longrightarrow \quad
0\le K\le v_p(a)+1. 
$$ 
\par
\noindent
Notice : the Hildebrand coefficients are, by definition, supported over the square-full numbers $n$: 
$$
n\enspace {\rm is}\enspace \hbox{\stampatello square-full} 
\enspace \definiz \enspace 
p|n \enspace \Rightarrow \enspace p^2|n. 
$$
\par
\noindent
It is a matter of convenience, then, to define $1$ as a square-full number; it's strange, since it is also square-free (both coming from $\omega(1)=0$, no prime-divisors!), but we'll do it, in order to simplify many details. In fact, $\Hildebrand(1)=F(1)$, of course, might be non-zero (thus confirming $1$ in the $\Hildebrand$ support, of square-full numbers). 
\par
\noindent
The vertical limit of $c_{p^{v_p(q)}}(a)$ (recall $c_q(a)$ is $q-$multiplicative, so may consider the vanishing of this factor), when combined to a square-full supported coefficient (with same prime $p|q$) $G(q)=G(p^{v_p(q)}\cdot q/p^{v_p(q)})$, gives: 
$$
2\le v_p(q)\le v_p(a)+1, 
\quad \forall p|q, 
$$
\par
\noindent
first inequality from assuming $G$ supported on square-full (and this $p$ dividing modulus $q$), second one to avoid $c_q(a)=0$ (from Ramanujan vertical limit). 
\par				
This, in turn, implies $v_p(a)\ge 1$, i.e. $p|a$ : we are proving that $p|q$ implies $p|a$ and so (assuming $G$, even not multiplicative, vanishing outside the square-full numbers) 
$$
\RG a = \sum_{q=1}^{\infty}G(q)c_q(a) = \sum_{{q=1}\atop {p|q\,\Rightarrow\,p|a}}^{\infty}G(q)c_q(a). 
$$
\par
\noindent
However, previous inequalities for $p-$adic valuations may be written 
$$
2\le v_p(q)\le v_p(a\rad a), 
\quad \forall p|q 
$$
\par
\noindent
which, if combined with: $p|q\,\Rightarrow\,p|a$, gives $q|a\rad a$ : 
$$
\RG a = \sum_{q|a\rad a}G(q)c_q(a), 
\quad \forall a\in\N, 
$$
\par
\noindent
but also, from our hypothesis $G$ is square-full supported, that $q=m\rad m$, with $m|a$, giving 
$$
\RG a = \sum_{m|a}G(m\rad m)c_{m\rad m}(a), 
\quad \forall a\in\N. 
$$
\par
\noindent
Thus getting Hildebrand's finite Ramanujan expansion. Notice that we started from a generic Ramanujan expansion representation for our $F$ and assuming ONLY its coefficient $G$, as for $\Hildebrand$, supported on square-full numbers, we proved it's actually a {\it finite Ramanujan expansion} abbrev. {\stampatello f.R.e.} hefereafter. 
\smallskip
\par
In case : $F\neq \0$ is multiplicative, we have (appearently) ANOTHER {\stampatello f.R.e.}, namely that of coefficient $G=G_F$, the canonical Ramanujan coefficient of our $F$. 
\par
Now, before proceeding further : the canonical Ramanujan coefficients, as may be easily checked, of course DO NOT DEPEND on $a$ (say, outer variable) neither their support does depend on $a$. We call in [CM], [C1] and [C3] an expansion, with this property on the coefficients, a {\it pure Ramanujan expansion}; there the first Author THOUGHT that this coefficient $\Hildebrand$ was not pure ! Thanks to the second Author (that scrutinized much better Hildebrand's result, through its [ScSp] reproduction) a conclusion, on the {\bf purity of Hildebrand's coefficient}, is now in order !
\par
(Also, see that neither Theorem 3 in [C3] nor its proof have problems; whereas, the only problem is that, {\bf even if we have a pure \& finite Ramanujan expansion}, a bound on its length isn't always possible ! In fact, a priori all the {\stampatello pure f.R.e} in the argument $a\in\N$ have a length depending on all non-zero $p-$adic valuations of $a$; so, {\bf the length may be}, and usually is, {\bf unbounded} with $a$, i.e., as $a\to\infty$.)
\par
In case : $F\neq \0$ is multiplicative, then we know it has an expansion with coefficient $G=G_F$, the canonical Ramanujan coefficient of our $F$ (and it's a {\stampatello f.R.e.}). Well, our $G_F$, too, like $\Hildebrand$, is supported on square-full numbers (because it's multiplicative and vanishes on primes, so the $q$ in its support have $v_p(q)\ge 2$, for non-vanishing $p-$adic valuations). Let's see how we got $G_F$. We started from $G_F(1)=1$, $G_F(p)=0$ on all primes $p$, so $\SGf=\1$ and this means : 
$$
\corsivoR_{G_F}(a) = \sum_{q=1}^{\infty}G_F(q)c_q(a) 
= \sum_{d|a\rad a}G_F(d)c_d(a)\cdot \sum_{(r,a)=1}G_F(r)\mu(r)
= \sum_{d|a\rad a}G_F(d)c_d(a), 
\quad \forall a\in\N, 
$$
\par
\noindent
from the $\corsivoR \corsivoS-$formula (Corollary 3.3). We may proceed as above, since $G_F(d)\neq 0$ implies $d$ square-full; hence, $p|d$ implies $2\le v_p(d)\le v_p(a)+1$, like above, which is equivalent to: $1\le v_p(d/\rad d)\le v_p(a)$, $\forall p|d$; in which, now, setting $q:=d/\rad d$, whence (from $d$ square-full) $\rad q=\rad d$ transforms $d|a\rad a$ into $q|a$ and $d=q\rad q$: 
$$
\corsivoR_{G_F}(a) = \sum_{q|a}G_F(q\rad q)c_{q\rad q}(a), 
\quad \forall a\in\N. 
$$
\par
\noindent
Thus, any {\it multiplicative and square-free vanishing} Ramanujan coefficient of our multiplicative $F\neq \0$ is square-full supported, whence it's Hildebrand's coefficient of this F : $\Hildebrand$. A kind of astonishing, $G_F=\Hildebrand$. 
\par				
However, not so strange, if we think about HOW we defined $G_F$ : recursively, on the prime-powers; but this, in turn, adding multiplicativity, is the same as building $\Hildebrand$ recursively on divisor sums ! 
\par
\noindent
Summarizing, we have proved the following Lemmas, for finite Ramanujan expansions. 
\medskip
\par
First: of course, for any $F$, $\Hildebrand$ is the unique Ramanujan coefficient, of this $F$, with square-full support.
\smallskip
\par
\noindent {\bf Lemma A.6.1.} ({\stampatello uniqueness for the Hildebrand coefficient})
\par
\noindent
{\it Let } $F:\N\to\C$ {\it be any arithmetic function and let } $G\in<F>$ {\it be square-full supported. Then $G=\Hildebrand$. 
In case $F\neq \0$ is multiplicative, we get: $G_F=\Hildebrand$. Furthermore, adopting the convention $G_{\0}\defineq \0$, we have that \enspace $F$ \enspace is {\stampatello multiplicative iff} \enspace $\Hildebrand$ \enspace is {\stampatello multiplicative}. 
}
\medskip
\par
\noindent 
Notice : if $F=\0$, then $G$ square-full supported implies however $G={\rm Hi}_{\0}=\0$, since its support is empty (however, contained in the square-full numbers set) !
\medskip
\par
Second: any $F\neq \0$ has a pure {\stampatello f.R.e.}, since $\Hildebrand$ is pure and gives a {\stampatello f.R.e.}! 
\smallskip
\par
\noindent {\bf Lemma A.6.2.} ({\stampatello Ramanujan clouds are all non-empty, like their pure and finite subsets})
\par
\noindent
{\it Let } $F:\N\to\C$ {\it be any arithmetic function. Then $\Hildebrand \in <F>$. In particular, all Ramanujan clouds are non-empty.  
In case $F\neq \0$, we have a {\stampatello pure and finite Ramanujan expansion} for $F$, namely that with Ramanujan coefficient $\Hildebrand$; in other words, $F\neq \0$ $\Rightarrow $ $<F>_{\hbox{\stampatello pure\& fin}}\neq \emptyset$. For $F=\0$, if we consider the trivial Ramanujan expansion with $\0$ Ramanujan coefficient a pure \& finite one, then $<\0>_{\hbox{\stampatello pure\& fin}}\neq \emptyset$. However, Ramanujan's $G_R(q)\defineq 1/q$ $[R]$ proves that the pure part of $\0$ cloud is, say, $<\0>_{\hbox{\stampatello pure}}\neq \emptyset$. 
}
\medskip
\par
\noindent
(Notation $<F>_{\hbox{\stampatello pure}}$ here is $<F>_{\ast}$ in [C3], where $<F>_{\ast}\cap <F>_{\#}$ is written $<F>_{\hbox{\stampatello pure\& fin}}$ here.) 

\bigskip

\par
This Lemma makes a natural question arise : are there Ramanujan expansions of $\0$ which are pure, finite and non-trivial (i.e., with coefficient $G\neq \0$) ? 

\vfill
\eject

\par				
\centerline{\stampatello Bibliography}

\bigskip

\item{[Ca]}  R. Carmichael, {\sl Expansions of arithmetical functions in infinite series}, Proc. London Math.
Society, 34 (1932), 1-26.
\smallskip
\item{[C1]} G. Coppola, {\sl A map of Ramanujan expansions}, ArXiV:1712.02970v2. (Second Version) 
\smallskip
\item{[C2]} G. Coppola, {\sl Finite and infinite Euler products of Ramanujan expansions}, ArXiV:1910.14640v2 (Second Version) 
\smallskip
\item{[C3]} G. Coppola, {\sl Recent results on Ramanujan expansions with applications to correlations}, Rend. Sem. Mat. Univ. Pol. Torino {\bf 78.1} (2020), 57--82. 
\smallskip
\item{[CG]} G. Coppola and L. Ghidelli, {\sl Multiplicative Ramanujan coefficients of null-function}, ArXiV:2005.14666v2 (Second Version)  
\smallskip
\item{[CMS]} G. Coppola, M. Ram Murty and B. Saha, {\sl Finite Ramanujan expansions and shifted convolution sums of arithmetical functions},  J. Number Theory {\bf 174} (2017), 78--92. 
\smallskip
\item{[CM]} G. Coppola and M. Ram Murty, {\sl Finite Ramanujan expansions and shifted convolution sums of arithmetical functions, II}, J. Number Theory {\bf 185} (2018), 16--47. 
\smallskip
\item{[D]} H. Davenport, {\sl Multiplicative Number Theory}, 3rd ed., GTM 74, Springer, New York, 2000. 
\smallskip
\item{[De]} H. Delange, {\sl On Ramanujan expansions of certain arithmetical functions}, Acta Arith., 31 (1976), 259--270.
\smallskip
\item{[GrKnP]} R.L. Graham, D.E. Knuth and O. Patashnik, {\sl Concrete mathematics: a foundation for computer science}, Computers in Physics, 3(5), 106-107.
\smallskip
\item{[H]} G.H. Hardy, {\sl Note on Ramanujan's trigonometrical function $c_q(n)$ and certain series of arithmetical functions}, Proc. Cambridge Phil. Soc., {\bf 20} (1921), 263--271. 
\smallskip
\item{[HL]} G.H. Hardy and J.E. Littlewood, {\sl Contributions to the theory of the Riemann zeta-function and the theory of the distribution of primes}, Acta Mathematica {\bf 41} (1916), 119--196.
\smallskip
\item{[Hau]} P. Haukkanen, {\sl Extensions of the class of multiplicative functions}, ArXiV:1308.6670v1
\smallskip
\item{[Hi]} A. Hildebrand, {\sl \"Uber die punktweise Konvergenz von Ramanujan-Entwicklungen zahlentheoretischer Funktionen} [[On the pointwise convergence of Ramanujan expansions of number-theoretic functions]], Acta Arith. {\bf 44} (1984), no. 2, 109--140. {\tt MR}$0774094$ 
\smallskip
\item{[H\"o]} O. H\"older, {\sl Zur Theorie der Kreisteilungsgleichung $K_m (x)=0$}, Prace Mat. Fiz. 43, 13–23 (1936).
\smallskip
\item{[Hoo]} C. Hooley, {\sl A note on square-free numbers in arithmetic progressions}, Bull. London Math. Soc. {\bf 7} (1975), 133-138.
\smallskip
\item{[K]} J.C. Kluyver, {\sl Some formulae concerning the integers less than $n$ and prime to $n$}, Proceedings of the Royal Netherlands Academy of Arts and Sciences (KNAW), 9(1):408--414, 1906. 
\smallskip
\item{[Lah]} D.B. Lahiri, {\it Hypo-multiplicative number-theoretic functions}, Aequationes Math. {\bf 9} (1973), 184--192. 
\smallskip
\item{[La]} M. Laporta, {\sl On Ramanujan expansions with multiplicative coefficients}, to appear on Indian J. of Pure and Applied Math. 
\smallskip
\item{[Lu1]} L.G. Lucht, {\sl Ramanujan expansions revisited}, Archiv der Mathematik {\bf 64.2} (1995), 121--128.
\smallskip
\item{[Lu2]} L.G. Lucht, {\sl A survey of Ramanujan expansions}, International Journal of Number Theory {\bf 6} (2010), 1785--1799.
\smallskip
\item{[M]} M. Ram Murty, {\sl Ramanujan series for arithmetical functions}, Hardy-Ramanujan J. {\bf 36} (2013), 21--33. Available online 
\smallskip
\item{[Pa]} F. Pappalardi, {\sl A survey on $k$-freeness}, Proceeding of the Conference in Analytic Number Theory in Honor of Prof. Subbarao at IM Sc. Chennai. 2003.
\smallskip
\item{[R]} S. Ramanujan, {\sl On certain trigonometrical sums and their application to the theory of numbers}, Transactions Cambr. Phil. Soc. {\bf 22} (1918), 259--276.
\smallskip
\item{[Rr]} O. Ramar\'e, {\sl Explicit estimates on the summatory functions of the M\"obius function with coprimality restrictions}, Acta Arith. {\bf 165} (2014), 1--10.
\smallskip
\item{[Re]} D. Rearick, {\sl Semi-multiplicative functions}, Duke Math. J. {\bf 33} (1966), 49--53. 
\smallskip			
\item{[ScSp]} W. Schwarz and J. Spilker, {\sl Arithmetical Functions}, Cambridge University Press, 1994.
\smallskip
\item{[T]} G. Tenenbaum, {\sl Introduction to Analytic and Probabilistic Number Theory}, Cambridge Studies in Advanced Mathematics, {46}, Cambridge University Press, 1995. 
\smallskip
\item{[Wa]} A. Walfisz, {\sl Weylsche Exponentialsummen in der neueren Zahlentheorie}, Mathematische Forschungsberichte, XV. VEB Deutscher Verlag der Wissenschaften, Berlin 1963.
\smallskip
\item{[Wi]} A. Wintner, {\sl Eratosthenian averages}, Waverly Press, Baltimore, MD, 1943. 

\bigskip
\bigskip
\bigskip

\par
\leftline{\tt Giovanni Coppola - Universit\`{a} degli Studi di Salerno (affiliation)}
\leftline{\tt Home address : Via Partenio 12 - 83100, Avellino (AV) - ITALY}
\leftline{\tt e-mail : giovanni.coppola@unina.it}
\leftline{\tt e-page : www.giovannicoppola.name}
\leftline{\tt e-site : www.researchgate.net}

\bigskip
\bigskip

\leftline{\tt Luca Ghidelli - Universit\`{a} degli Studi di Genova (affiliation)}
\leftline{\tt Address: Via Dodecaneso 35 - CAP 16146, Genova (GE) - ITALY}
\leftline{\tt e-mail : luca.ghidelli@uottawa.ca}
\leftline{\tt e-page : math-lucaghidelli.site}
\leftline{\tt e-site : www.researchgate.net}

\bye